\documentclass[10pt]{amsart}
\newcommand\version{public}  
\newcommand\finalized{yes}
\newcommand\theoremnumbering{section}  
\usepackage{fullpage,setspace}
\usepackage{mathrsfs} 
\newcommand\choosefont[1]{\usepackage{#1}}
\usepackage[bottom]{footmisc}  
\usepackage{enumerate}   
\usepackage{url,cite,fancyheadings}
\usepackage{enumerate}
\usepackage{asymptote}
\usepackage{amsmath,amssymb,amscd,mathrsfs,latexsym}
\usepackage{mathdots}

\usepackage{ifthen}     

\newcommand\pubpri[2]{%
\ifthenelse{\equal{\version}{public}}%
{{#1}}%
{\ifthenelse{\equal{\finalized}{no}}{\marginpar{\scshape\small Pubpri Alert}{#2}}{#2}}{}}
\newcommand\pubprinoalert[2]{%
\ifthenelse{\equal{\version}{public}}%
{{#1}}%
{#2}}

\newcommand\ignore[1]{}

\usepackage{color}
\providecommand\wantcolor{yes}   %
\definecolor{backgroundyellow}{cmyk}{.2,.1,.8,.2}
\definecolor{backgroundblue}{rgb}{0,0,1}
\definecolor{backgroundred}{rgb}{1,0,0}
\definecolor{backgroundmagenta}{cmyk}{0,1,0,0}
\ifthenelse{\equal{\wantcolor}{no}}{\renewcommand\color[1]{}}{}

\newcommand\mysection{\section}

\newcommand\mysubsection{\subsection}

\newcommand\mysubsubsection[1]{%
		\subsubsection{\sffamily\upshape\mdseries #1}}

\newcommand\mysss{\mysubsubsection}

\usepackage{chngcntr}
\counterwithin{figure}{section}

\usepackage{url,hyperref}  
\usepackage{hyperref}
\hypersetup{colorlinks=true,
linkcolor=blue,
filecolor=magenta,
urlcolor=cyan,}
\providecommand{\theoremnumbering}{document}

\ifthenelse{\equal{\theoremnumbering}{section}}{\newtheorem{annotation}{Annotation}[section]}%
{\ifthenelse{\equal{\theoremnumbering}{subsection}}{\newtheorem{annotation}{Annotation}[subsection]}{\newtheorem{annotation}{Annotation}}}
\newtheorem{theorem}[annotation]{
		Theorem}
\newtheorem{lemma}[annotation]{
		Lemma}
\newtheorem{definition}[annotation]{
		Definition}
\newtheorem{corollary}[annotation]{
		Corollary}
\newtheorem{proposition}[annotation]{
		Proposition}
\newtheorem{conjecture}[annotation]{
		Conjecture}
\newtheorem{example}[annotation]{
		Example}
\newcommand\bexample{\begin{example}\begin{rm}}
\newcommand\eexample{\end{rm}\hfill$\Box$\end{example}}
\newtheorem{examplenobox}[annotation]{
		Example}
\newcommand\bexamplenobox{\begin{examplenobox}\begin{rm}}
\newcommand\eexamplenobox{\end{rm}\end{examplenobox}}
\newtheorem{exercise}[annotation]{
		Exercise}
\newcommand\bexercise{\begin{exercise}\begin{rm}}
\newcommand\eexercise{\end{rm}\end{exercise}}
\newtheorem{notation}[annotation]{
		Notation}
\newcommand\bnotation{\begin{notation}\begin{rm}}
\newcommand\enotation{\end{rm}\end{notation}}

\newtheorem{remark}[annotation]{
		Remark}
\newcommand\bremark{\begin{remark}
\begin{upshape}}
\newcommand\eremark{\end{upshape}
\end{remark}}
\newenvironment{remark*}{%
\par\noindent{\scshape 
  Remark: }\begin{rm}}{\hfill\end{rm}\newline} 
\newcommand\bremarkstar{\begin{remark*}}
\newcommand\eremarkstar{\end{remark*}}
\newcommand\bdefn{\begin{definition}
\begin{upshape}}
\newcommand\edefn{\end{upshape}
\end{definition}}
\newtheorem{caveat}[annotation]{
		Caveat}
\newcommand\bcaveat{\begin{caveat}
\begin{upshape}}
\newcommand\ecaveat{\end{upshape}
\end{caveat}}
\newenvironment{caveatstar}{
\par\noindent{\scshape\bfseries
  Caveat: }\begin{rm}}{\end{rm}\newline} 
\newcommand\bcaveatstar{\begin{caveatstar}}
\newcommand\ecaveatstar{\end{caveatstar}}
\newenvironment{myproof}{%
\par\noindent{\scshape 
  Proof: }\begin{rm}}{\hfill$\Box$\end{rm}\newline} 
\newcommand\bmyproof{\begin{myproof}}
\newcommand\emyproof{\end{myproof}}
\newenvironment{myproofnobox}{%
\par\noindent{\scshape Proof: }\begin{rm}}{\end{rm}\hfill\newline}
\newcommand\bmyproofnobox{\begin{myproofnobox}}
\newcommand\emyproofnobox{\end{myproofnobox}}
\newenvironment{solution}{%
\par\noindent{\scshape Solution: }\begin{rm}}{\hfill$\Box$\end{rm}\newline}
\newenvironment{solutionnobox}{%
\par\noindent{\scshape Solution: }\begin{rm}}{\end{rm}}
\newcommand\bsolution{\begin{solution}\begin{rm}}
\newcommand\esolution{\end{rm}\end{solution}}
\newcommand\bsolutionnobox{\begin{solutionnobox}\begin{rm}}
\newcommand\esolutionnobox{\end{rm}\end{solutionnobox}}
\newcommand\bthm{\begin{theorem}}
\newcommand\ethm{\end{theorem}}
\newcommand\bcor{\begin{corollary}}
\newcommand\ecor{\end{corollary}}
\newcommand\blemma{\begin{lemma}}
\newcommand\elemma{\end{lemma}}
\newcommand\bprop{\begin{proposition}}
\newcommand\eprop{\end{proposition}}
\newcommand\beqn{\begin{equation}}
\newcommand\eeqn{\end{equation}}
\newcommand\beqnstar{\begin{equation*}}
\newcommand\eeqnstar{\end{equation*}}
\usepackage{amsthm}

\newcommand\mtitle[1]%
{\marginpar{\sffamily\raggedright\upshape\small #1}}

\providecommand\finalized{no}
\ifthenelse{\equal{\finalized}{yes}}{}{\marginparwidth=2cm}
\ifthenelse{\equal{\finalized}{yes}}%
{\newcommand\mylabel[1]{\label{#1}}}%
{\newcommand\mylabel[1]{\label{#1}\marginpar{[{\ttfamily\upshape\tiny #1}]}}}
\ifthenelse{\equal{\finalized}{yes}}%
{\newcommand\checked[1]{}}%
{\newcommand\checked[1]{\marginpar{[{\ttfamily\upshape\tiny CHECKED: #1}]}}}
\ifthenelse{\equal{\finalized}{yes}}%
{\newcommand\spellchecked[1]{}}%
{\newcommand\spellchecked[1]{\marginpar{[{\ttfamily\upshape\tiny SPELLCHECKED: #1}]}}}

\providecommand\version{public}   
\ifthenelse{\equal{\version}{public}}%
{\newcommand\mcomment[1]{}}%
{\newcommand\mcomment[1]{\marginpar{{\raggedright\sffamily\upshape\small
\begin{spacing}{0.75} #1\end{spacing}}}}}
\ifthenelse{\equal{\version}{public}}%
{\newcommand\fcomment[1]{}}%
{\newcommand\fcomment[1]{\footnote{#1}}}
\ifthenelse{\equal{\version}{public}}%
{\newcommand\comment[1]{}}%
{\newcommand\comment[1]{{\small #1}}}

\usepackage{mathrsfs,graphicx,amsthm,amsfonts,amsxtra,xspace,array,cite,epic,float,ifpdf}
\usepackage{youngtab}
\choosefont{newcent}    

\title[on Chari-Loktev bases]{On Chari-Loktev bases for local Weyl
modules in type $A$}
\author{K.~N.~Raghavan}
	\address{The Institute of Mathematical Sciences, HBNI, 
                 C.~I.~T.~Campus, 
		Chennai 600\,113, INDIA}
	\email{knr@imsc.res.in}
\author{B.~Ravinder}
	\address{Tata Institute of Fundamental Research, 
                Homi Bhabha Road, 
		Mumbai 400\,005, INDIA}
	\email{ravinder@math.tifr.res.in}

\author{Sankaran Viswanath}
	\address{The Institute of Mathematical Sciences, HBNI, 
                C.~I.~T.~Campus, 
		Chennai 600\,113, INDIA}
	\email{svis@imsc.res.in}
\thanks{The second named author acknowledges support from CSIR under the SPM 
  doctoral fellowship scheme. The first and third named 
authors acknowledge support from DAE under a XII plan project}
\subjclass[2010]{17B67 (05E10)}

\newcommand\st{\,|\,}

\newcommand\lieg{\mathfrak{g}}
\newcommand\lieh{\mathfrak{h}}
\newcommand\lieb{\mathfrak{b}}

\newcommand\sltwo{\mathfrak{sl}_2}

\newcommand\tensor{\otimes}

\newcommand\aseq{\underline{a}}
\newcommand\bseq{\underline{b}}\newcommand\cseq{\underline{c}}
\newcommand\lseq{\underline{\lambda}}
\newcommand\kseq{\underline{\kappa}}
\newcommand\lpseq{\underline{\lambda}'}
\newcommand\kpseq{\underline{\kappa'}}
\newcommand\lseqlong{\lambda_1\geq\ldots\geq\lambda_n}

\newcommand\mseq{\underline{\mu}}

\newcommand\tseq{\underline{\theta}}

\newcommand\murpone{\mu_{r+1}}
\newcommand\mutwo{\mu_2}

\newcommand\lj{\lambda^j}
\newcommand\ljpone{\lambda^{j+1}}

\newcommand\mynum{\mathfrak{m}}

\newcommand\pijiseq{\underline{\pi(j)}^i}
\newcommand\pijoneseq{\underline{\pi(j)}^1}
\newcommand\pijseq{\underline{\pi(j)}}
\newcommand\pirseq{\underline{\pi(r)}}
\newcommand\pioneseq{\underline{\pi(1)}}
\newcommand\pitwoseq{\underline{\pi(2)}}
\newcommand\ljseq{\lseq^j}
\newcommand\ljponeseq{\lseq^{j+1}}

\DeclareMathOperator\depth{\textup{depth}}

\newcommand\intpart{\mathcal{M}}
\newcommand\npattern{\mathscr{N}}
\newcommand\nplrpone{\npattern_{\lseqrpone}}
\newcommand\nlmd{\npattern_{\lseq,\mseq}(d)}
\newcommand\nlmkd{\npattern_{\lseq,\mseq}^k(d)}
\newcommand\plmkd{\popset_{\lseq,\mseq}^k(d)}

\newcommand\plmkpjsd{\popset_{\lseq,\mseq}^{k+j}[d]}
\newcommand\plmksd{\popset_{\lseq,\mseq}^k[d]}

\newcommand\comp{\complement}

\newcommand\popset{\mathbb{P}}
\newcommand\one{\underline{\mathbf{1}}}

\newcommand\pkl{\popset^k_{\lseq}}
\newcommand\pil{\popset^i_{\lseq}}

\newcommand\pipjl{\popset^{i+j}_{\lseq}}
\newcommand\pnotl{\popset^0_{\lseq}}
\newcommand\ptwol{\popset^2_{\lseq}}
\newcommand\ponel{\popset^1_{\lseq}}
\newcommand\popnot{\pop_0}
\newcommand\vcl[1]{v_{#1}}
\newcommand\rhocl[1]{\rho_{#1}}
\newcommand\popnotk{\pop_0^k}

\newcommand\least{b}

\newcommand\crpone{\mathbb{C}^{r+1}}

\newcommand\philm{\Phi_{\lseq,\mseq}}
\newcommand\shift{\mathcal{S}}
\newcommand\shiftk{\shift^k}

\newcommand\shiftjpk{\shift^{j+k}}
\newcommand\shiftj{\shift^{j}}
\newcommand\phishift{\Phi_{\llseq,(\mu_1+k,\ldots,\mu_{r+1}+k)}}

\newcommand\pic{\pi^{\textup{c}}}
\newcommand\picseq{\piseq^{\textup{c}}}

\newcommand\vlseq{V(\lseq)}
\newcommand\vmseq{V(\mseq)}
\newcommand\piseq{\underline{\pi}}
\newcommand\piseqj{{\piseq}^j}
\newcommand\piseqone{{\underline{\pi}^1}}\newcommand\piseqtwo{{\underline{\pi}}^2}

\newcommand\ta{\underline{t}^{\underline{a}}}

\newcommand\pat{\mathbb{\pi}}
\newcommand\gtp{\mathfrak{P}}
\newcommand\gtpint{\mathfrak{P}^{\textup{int}}}

\newcommand\gtpintlseq{\gtpint_{\lseq}}
\DeclareMathOperator{\weight}{wt}
\newcommand\tul{\underline{t}}

\newcommand\xul{\underline{x}}
\newcommand\xntuple{(x_1,\ldots,x_n)}
\newcommand\xuldown{\xul^\downarrow}
\newcommand\xntupledown{(x_1^\downarrow,\ldots,x_n^\downarrow)}
\newcommand\yul{\underline{y}}
\newcommand\ydown{y^\downarrow}
\newcommand\xdown{x^\downarrow}
\newcommand\Rn{\mathbb{R}^n}

\newcommand\wmajgeq{\succcurlyeq_{\textup{wm}}}
\newcommand\maj{\preccurlyeq_{\textup{m}}}
\newcommand\majleq{\maj}
\newcommand\majgeq{\succcurlyeq_{\textup{m}}}
\newcommand\interlace{\gtrless}

\newcommand\herm{\mathcal{H}}
\newcommand\hermlseq{\herm_{\lseq}}

\newcommand\symm{\mathcal{S}}
\newcommand\symmlseq{\symm_{\lseq}}

\newcommand\specdown{\lseq^\downarrow}

\newcommand\Dlseq{D_{\lseq}}
\newcommand\nwsm[1]{\underline{#1}|}

\DeclareMathOperator\mtopat{\gtp}
\DeclareMathOperator\optraparea{\square}
\DeclareMathOperator\opproptraparea{\square^\textup{prop}}
\DeclareMathOperator\oparea{\bigtriangleup}

\newcommand\traparea[1]{\optraparea(#1)}
\newcommand\proptrap[1]{\opproptraparea(#1)}
\newcommand\area[1]{\oparea(#1)}

\newcommand\linear{T}
\newcommand\linearn{\linear(n)}
\newcommand\linearnmone{\linear(n-1)}

\newcommand\lseqn{{\lseq}^n}
\newcommand\lambdan{{\lambda}^n}
\newcommand\lseqnmone{\lseq^{n-1}}
\newcommand\lambdanmone{\lambda^{n-1}}
\newcommand\ltseq{{\tilde{\lseq}}^n}
\newcommand\lambdatn{{\tilde{\lambda}}^n}

\newcommand\kseqnmone{\kseq^{n-1}}
\newcommand\knmone{\kappa^{n-1}}

\newcommand\lseqj{\lseq^j}
\newcommand\lseqjmone{\lseq^{j-1}}
\newcommand\kseqj{\kseq^j}
\newcommand\kappan{\kappa^n}
\newcommand\kseqn{\kseq^n}

\newcommand\knot{{j_0}}
\newcommand\lnmone{\lambda^{n-1}}
\newcommand\lambdaj{\lambda^j}
\newcommand\ljmone{\lambda^{j-1}}

\newcommand\winterlace{{\interlace}_w}

\newcommand\etatwo{\eta'}
\newcommand\etaseq{\underline{\eta}}
\newcommand\etanmone{\etaseq^{n-1}}
\newcommand\etatwoseq{\underline{\etatwo}}

\newcommand\lseqjpone{\lseq^{j+1}}
\newcommand\kseqjpone{\kseq^{j+1}}

\newcommand\slrpone{\mathfrak{sl}_{r+1}}
\newcommand\currlieg{\lieg[t]}
\newcommand\curralg\currlieg

\newcommand\wm{W}
\newcommand\wml{\wm(\lambda)}
\newcommand\nplus{\mathfrak{n}^+}
\newcommand\nminus{\mathfrak{n}^-}
\newcommand\ghat{\hat{\lieg}}
\newcommand\hhat{\hat{\lieh}}
\newcommand\bhat{\hat{\lieb}}

\newcommand\germ\mathfrak

\newcommand\wmlmM{W(\lambda)_\mu[M]}
\newcommand\dem{V_w(\Lambda)}
\newcommand\wmlseq{W(\lseq)}

\newcommand\hwel{u}
\newcommand\pop{\mathfrak{P}}

\newcommand\bijxi{\Xi}
\newcommand\bijxieta{\bijxi_{\etaseq}}
\newcommand\bijxipeta{\bijxi'_{\etaseq}}
\newcommand\mustar{\mu^\star}
\newcommand\pistar{\pi^\star}
\newcommand\piseqstar{\piseq^\star}
\newcommand\pistarseq{\piseqstar}
\newcommand\astar{a^\star}
\newcommand\bstar{b^\star}
\newcommand\astarseq{\aseq^\star}
\newcommand\bstarseq{\bseq^\star}
\newcommand\lrpone{\lambda^{r+1}}
\newcommand\lseqrpone{\lseq^{r+1}}

\newcommand\lrponeseq{\lseqrpone}
\newcommand\loneseq{\lseq^1}

\newcommand\bijxilrpone{\bijxi_{\lrponeseq}}
\newcommand\bijxiljpone{\bijxi_{\ljponeseq}}
\newcommand\bijxipljpone{\bijxi'_{\ljponeseq}}
\newcommand\bijxiplrpone{\bijxi'_{\lrponeseq}}

\newcommand\lrseq{\lseq^r}
\newcommand\lr{\lambda^r}

\newcommand\pattern{\mathcal{P}}
\newcommand\leta{\tilde{\eta}}
\newcommand\letatwo{\tilde{\eta}'}
\newcommand\letaseq{\tilde{\etaseq}}
\newcommand\letatwoseq{\underline{\tilde{\eta}}'}

\newcommand\rpartd{\mathscr{P}_r(d)}

\newcommand\llr{\tilde{\lambda}^r}
\newcommand\llseq{\underline{\tilde{\lambda}}}
\newcommand\llrseq{\tilde{\underline{\lambda}}^r}
\newcommand\llrponeseq{\tilde{\lseq}^{r+1}}
\newcommand\patone{\pattern_1}
\newcommand\lpattern{\tilde{\pattern}}
\newcommand\barpop{\overline{\pop}}
\newcommand\pijciseq{\underline{\pi^{\textup{c}}(j)}^i}

\newcommand\beq{\begin{equation}}
\newcommand\eeq{\end{equation}}

\newcommand{\be}{\begin{enumerate}}
\newcommand{\ee}{\end{enumerate}}

\newcommand{\integers}{\mathbb{Z}}

\renewcommand\omega{\varpi}

\newcommand{\bdf}{\mathbf{F}}
\newcommand{\bds}{\mathbf{s}}

\newcommand{\clpop}{\mathfrak{C}}

\begin{document}
\begin{abstract} 
This paper is a study of the bases introduced by Chari-Loktev in~\cite{cladv2006} for local Weyl modules of the
current algebra associated to a special linear Lie algebra.
Partition overlaid patterns, POPs for short---whose introduction is
one of the aims of this paper---form  convenient parametrizing sets of these bases.
They play a role analogous to that played by (Gelfand-Tsetlin) patterns in the
representation theory of the special linear Lie algebra.

The notion of a POP leads naturally to the
notion of area of a pattern.     We observe that there is a unique
pattern of maximal area among all those with a given bounding sequence
and given weight.   We give a combinatorial proof of this and discuss
its representation theoretic relevance.    

We then state a conjecture about the ``stability'', i.e.,  compatibility in the long range, of Chari-Loktev bases with respect to
inclusions of local Weyl modules.   In order to state the conjecture,
we establish a certain bijection between colored partitions and POPs,
which may be of interest in itself.    The stability conjecture has
been proved in~\cite{rrv:stab} in the rank one case.    
%
\end{abstract}
\keywords{Gelfand-Tsetlin pattern, Chari-Loktev basis, Demazure
  module, partition overlaid pattern, 
  POP, area of a pattern, stability of Chari-Loktev bases}

\maketitle
\mysection{Introduction}\mylabel{s:introduction}
\noindent
This paper is a study of the bases introduced by
Chari-Loktev~\cite{cladv2006} for local Weyl modules of current algebras of type $A$.  
Let $\lieg=\slrpone$, the simple Lie algebra of traceless complex $(r+1)\times(r+1)$ matrices.
Let $\currlieg=\lieg\tensor\mathbb{C}[t]$ denote the corresponding current algebra,  namely,  the extension by scalars of~$\lieg$ to the polynomial ring $\mathbb{C}[t]$.    We think of $\currlieg$ as a Lie algebra over the complex numbers, graded by $t$, and are interested in the representation theory of its graded finite dimensional 
modules.


Fix subalgebras of $\lieg$ as follows:  a Cartan subalgebra $\lieh$ and a Borel subalgebra $\lieh\oplus\nplus$ 
containing~$\lieh$. 
We may, without loss of generality, take these to be the diagonal and upper triangular subalgebras respectively.    
  It is easy to see that a graded (non-zero) finite dimensional
  $\currlieg$-module admits a homogeneous {\em highest weight
  vector\/} $u\neq0$, in the sense that there exists an element
  $\lambda$ of $\lieh^\star$ such that:
  \begin{equation}\label{e:intro1}  
  \quad (\nplus\tensor\mathbb{C}[t])\,\hwel=0,\quad\quad (\lieh\tensor
  t\mathbb{C}[t])\,\hwel=0,\quad\quad
  H\,\hwel=\langle\lambda,H\rangle\hwel\quad\textup{ for all $H\in\lieh$} 
%
\end{equation}

    The weight $\lambda$ in (\ref{e:intro1}) is forced to be dominant integral (because it is a highest weight for the restriction to~$\lieg$).  Furthermore,  for any simple root $\alpha$,   for $Y_\alpha$ in the root space $\lieg_{-\alpha}$, we must have:
\begin{equation}\label{e:intro2}
Y_\alpha^{\langle \lambda, \alpha^\vee\rangle+1}u=0
\end{equation}
where~$\alpha^\vee$ denotes the coroot corresponding to~$\alpha$.
Thus it is natural to define a module generated by a homogeneous element $u$ subject to the relations~(\ref{e:intro1}), (\ref{e:intro2}).      This is the {\em local Weyl module\/} associated to the dominant integral weight $\lambda$ (defined in~\cite{cprt2001}).  It turns out to be itself finite dimensional.     We denote it by~$\wml$.


Let a dominant integral weight $\lambda$ be fixed for the rest of this section.
Chari-Loktev~\cite{cladv2006}   construct a basis of homogeneous weight vectors for the local Weyl module $\wml$,  which is analogous in many ways to the well known Gelfand-Tsetlin basis for the irreducible $\lieg$-module $V(\lambda)$ with highest weight~$\lambda$.

One of our aims in the present paper 
is to clarify the parametrizing set of this Chari-Loktev basis.     Towards this end we introduce the notion of a {\em partition overlaid pattern\/} or {\em POP\/} (see~\S\ref{ss:pop}).
Dominant integral weights for $\lieg=\slrpone$ may be identified with
non-increasing sequences of non-negative integers of length $r+1$ with
the last element of the sequence being $0$.   Gelfand-Tsetlin patterns
(or just patterns---see~\S\ref{s:notation} for the precise definition)
with bounding sequence (corresponding to)~$\lambda$ parametrize the
Gelfand-Tsetlin basis for~$V(\lambda)$.     Analogously, as we show in~\S\ref{ss:clbase},
POPs with bounding sequence~$\lambda$ parametrize the Chari-Loktev basis for $\wml$.
The weight of the underlying pattern of a POP equals the $\lieh$-weight of the corresponding basis element and the number of boxes in the partition overlay determines the grade.   

The notion of a POP leads naturally to the notion of the {\em area\/} of a pattern (see~\S\ref{ss:area}). For a weight $\mu$ of~$V(\lambda)$,  it turns out that the piece of highest grade in the $\mu$-weight space of the local Weyl module~$\wml$  is one dimensional (the $\lieh$-weights of $\wml$ are precisely those of~$V(\lambda)$).    We give a representation theoretic proof of this fact in~\S\ref{s:relevance}.    This suggests---even proves, albeit circuitously---that there must be a unique pattern of highest area among all those with bounding sequence~$\lambda$ and weight~$\mu$.   We give a direct, elementary, and purely combinatorial proof of this in~\S\ref{s:theorem}.      

In~\S\ref{s:stabnew},  we state a conjecture about the ``stability'' of the Chari-Loktev basis.   To describe what is meant by stability,  let $\theta$ be the highest root of~$\lieg$.     We then have natural inclusions of local Weyl modules:
\[ \wml\hookrightarrow W(\lambda+\theta)\hookrightarrow W(\lambda+2\theta)\hookrightarrow\ldots \]
We may ask if there are corresponding natural inclusions of indexing sets of the Chari-Loktev bases of these modules.    To prove that this is indeed the case,  we first establish a combinatorial bijection (Theorem~\ref{t:bijection}),  which is a generalization of the construction of Durfee squares.
Loosely stated,  the bijection identifies colored partitions of an integer with a certain set of POPs.

Once the inclusion of indexing sets is established,   we may ask if
the Chari-Loktev bases respect the inclusions.   We believe that they
do have this {\em stability\/} property. 
More concretely, we make a precise conjecture to this effect (see \S\ref{s:stabnew}).       In earlier work~\cite{rrv:stab},   we have proved the conjecture in the case $r=1$ (namely for $\sltwo$).\footnote{The general case has since been proved by one of us~\cite{br:stab}.}


This paper is so arranged that the purely combinatorial sections
(\S\ref{s:theorem}, \S\ref{s:bijection}) can
be read independently of the rest.

\medskip
\noindent
{\em Acknowledgments:}
The authors thank Vyjayanthi Chari for many helpful discussions. Some of the initial explorations leading up to this work were carried out using {\tt Sage} \cite{sagemath}; it is our pleasure to thank the SageMath development team.

\mysection{Notation}\mylabel{s:notation}
\noindent 
This section establishes notation and terminology.   The notion of a
partition overlaid pattern, or POP, is introduced in~\S\ref{ss:pop}.    POPs parametrize Chari-Loktev bases for
local Weyl modules of~$\slrpone[t]$ (just as patterns parametrize the Gelfand-Tsetlin
bases for irreducible representations of~$\slrpone$).     The notion of area
(triangular and trapezoidal) of a pattern introduced
in~\S\ref{ss:area} plays an important role in what follows.
\mysubsection{Interlacing condition on sequences}\mylabel{ss:interlace}
Let $\lseq$ be a sequence $\lambda_1$, \ldots, $\lambda_n$ of real numbers.
\begin{itemize}
\item $\lseq$ is {\em (non-negative) integral\/} if the $\lambda_j$
are all (non-negative)
integers.    
\item
The number $n$ of elements in $\lseq$ is called its {\em length\/}.
\item
Let  $\lseq$: $\lseqlong$ and $\mseq$: $\mu_1\geq\ldots\geq\mu_{n-1}$ be non-increasing sequences of lengths~$n$ and~$n-1$ respectively.
We say that they {\em interlace\/} and write $\lseq\interlace\mseq$ if
\begin{equation}
{\lambda_1\geq\mu_1\geq\lambda_2}, \quad\quad
{\lambda_2\geq\mu_2\geq\lambda_3}, \quad\quad 
\ldots,\quad\quad\textup{and}\quad\quad
{\lambda_{n-1}\geq\mu_{n-1}\geq\lambda_n}.
\end{equation}
The interlacing condition may be remembered easily if $\lseq$ and
$\mseq$ are arranged like so:
\begin{equation}
\begin{array}{ccccccccccccccc}
&\mu_1 & &\mu_2 & & \cdots & & \cdots & & \cdots & & \mu_{n-2} & & \mu_{n-1}\\
\lambda_1 & &\lambda_2 & & \lambda_3 & & \cdots & & \cdots & &\lambda_{n-2}
& & \lambda_{n-1} & & \lambda_{n}
\end{array}
\end{equation}
If we now imagine $\geq$ relations among numbers as we move in the
north-easterly or the south-easterly direction,  that is precisely the
condition for interlacing.
\item   
For $\lseq$ and $\mseq$ as in the previous item,   we will feel free to use several alternative
  expressions to express the condition that they interlace:
\begin{quote}
  $\mseq\interlace\lseq$; \quad   $\mseq$ interlaces $\lseq$;\quad
  $\lseq$ interlaces $\mseq$; \quad $\lseq$ and $\mseq$ are
  interlaced;  \quad etc.
\end{quote}
\end{itemize}
\mysubsubsection{Weak interlacing}\mylabel{sss:wilace}
Let $\lseq: \lambda_1\geq\cdots\geq\lambda_n$ be a non-increasing
sequence of real numbers and $\mseq=(\mu_1,\ldots,\mu_{n-1})$ be an
element of $\mathbb{R}^{n-1}$.    We say that  $\lseq$ {\em weakly
  interlaces\/} $\mseq$ and write $\lseq\winterlace\mseq$ if for every $j$, $1\leq j\leq n-1$, and every
sequence $1\leq i_1<\ldots<i_j\leq n-1$,   we have:
\begin{equation}\label{e:winterlace}
\lambda_1+\cdots+\lambda_j\geq \mu_{i_1}+\cdots+\mu_{i_j}\geq\lambda_{n-j+1}+\cdots+\lambda_n
\end{equation}
It is evident that, for non-increasing sequences $\lseq$ and $\mseq$
of lengths $n$ and $n-1$,   
if $\lseq$ interlaces $\mseq$ then $\lseq$ weakly interlaces $\mseq$.
\mysubsection{Gelfand-Tsetlin patterns}\mylabel{ss:gtpatterns}
A {\em Gelfand-Tsetlin pattern\/},  or {\em GT pattern\/}, or just
{\em pattern\/} is a finite sequence of interlacing sequences.   More
precisely,  a pattern consists of a sequence $\lseq^1$, \ldots,
$\lseq^n$ of sequences such that 
\begin{itemize}
\item $\lseq^j$ is of length $j$  for all $1\leq j\leq n$;
\item $\lseq^{j}$ interlaces $\lseq^{j+1}$ for $1\leq j\leq n-1$:
that is, $\lseq^1\interlace\lseq^{2}\interlace\cdots\interlace\lseq^{n-1}\interlace\lseq^{n}$.
 \end{itemize}
The sequences in a pattern are arranged one below the other,  in a staggered fashion.    For example, the pattern consisting of the sequences $5$; $7$, $4$; and $7$, $5$, $3$ is written:
\begin{equation}\label{e:patex}
\begin{array}{ccccc}
&&5&&\\
&7&&4\\
  7&&5&&3\\
\end{array}
\end{equation}

The last sequence of a pattern is its
 {\em bounding sequence}.    For instance, the bounding sequence  of
 the pattern (\ref{e:patex}) is $7$, $5$, $3$.       When we speak of
 a pattern $\lseq^1$, \ldots, $\lseq^n$,   it is often convenient to
 let $\lseq^0$ denote the empty sequence.


\mysubsubsection{Integral patterns}\mylabel{sss:intpat}
A pattern is {\em integral\/} if all its entries are integers.   
\mysubsubsection{Rows of a pattern}\mylabel{sss:rows}
Let $\pattern$ be the pattern $\lseq^1$, \ldots, $\lseq^n$.  The entries
of $\lseq^k$   are sometimes referred to as the {\em entries in row~$k$\/}
of~$\pattern$.  The lone entry in the first row of the pattern (\ref{e:patex}) 
is $5$; the entries in its second row are $7$, $4$;  and those in third are
$7$, $5$, and~$3$.  



\mysubsubsection{The weight of a pattern}\mylabel{sss:weight}  
The {\em weight\/} of a pattern with bounding sequence of length~$n$ is
the $n$-tuple $(a_1,\ldots, a_n)$, 
where $a_j$ is the difference of the sum of the entries in row $j$ and
the sum of the entries in row~$j-1$.   It is understood that the sum
of the entries in the zeroth row is zero.  The weight of the 
pattern in (\ref{e:patex}), for instance,  is $(5,6,4)$.      

\mysubsection{Trapezoidal Area and (Triangular) Area of a pattern}\mylabel{ss:area}
Let $\lseq$ and $\mseq$ be two non-increasing sequences of lengths $n$ and~$n-1$
respectively that are interlaced.      The {\em triangular 
  area\/} or just {\em area\/} of the pair $(\lseq,\mseq)$ is defined by:
\begin{equation}\label{e:defarea}
\area{\lseq,\mseq}:=\sum_{i=1}^{n-1}\ (\lambda_i-\mu_i)(\mu_i-\lambda_{i+1}) 
\end{equation}
And the {\em trapezoidal 
  area\/} of the pair $(\lseq,\mseq)$ is defined by:
\begin{equation}\label{e:deftraparea}
\traparea{\lseq,\mseq}:=\sum_{1\leq i\leq j\leq 
  n-1}(\lambda_i-\mu_i)(\mu_j-\lambda_{j+1})   
\end{equation}
The above definitions make sense even when $n=1$:
$\mseq$ is empty and both areas vanish (since they are
empty sums).

The {\em triangular area\/}  or just {\em area\/} of a pattern $\pattern$ with rows 
$\lseq^1$, \ldots, $\lseq^n$ is defined by:
\begin{equation}\label{e:areadef:pat}
\area{\pattern}:=\area{\lseq^n,\lseq^{n-1}}+\area{\lseq^{n-1},\lseq^{n-2}}+\cdots+\area{\lseq^2,\lseq^1}+\area{\lseq^1,\lseq^0}
\end{equation}
Its {\em trapezoidal area\/} is  defined by:
\begin{equation}\label{e:trapareadef:pat}
\traparea{\pattern}:=\traparea{\lseq^n,\lseq^{n-1}}+\traparea{\lseq^{n-2},\lseq^{n-2}}+\cdots+\traparea{\lseq^2,\lseq^1}+\traparea{\lseq^1,\lseq^0}
\end{equation}
Observe that $\traparea{\pattern}\geq\area{\pattern}$ and that both
areas are zero if $\pattern$ has a single row.


\mysubsection{Majorization}\mylabel{ss:major}
For an element $\xul=\xntuple$ in~$\mathbb{R}^n$,  let
$\xuldown=\xntupledown$  be the vector whose co-ordinates are obtained
by rearranging the $x_j$ in non-increasing order.        For elements
$\xul$ and $\yul$ in~$\Rn$,   we say that $\xul$  {\em weakly             
  majorizes\/}  $\yul$ and write $\xul\wmajgeq\yul$ if                  
\begin{equation}\label{e:d:maj}
\xdown_1+\cdots+\xdown_k\geq \ydown_1+\cdots+\ydown_k
\quad\quad\textup{for 
  all $1\leq k\leq n$}
\end{equation}
The right hand side in the above equation is evidently the largest
possible value of $\sum_{i=1}^ky_{j_i}$ over all sequences $1\leq
j_1<\ldots<j_k\leq n$.    Thus (\ref{e:d:maj}) is equivalent to the a
priori stronger condition:
\begin{equation}\label{e:d:maj:cor}
\xdown_1+\cdots+\xdown_k
\geq
y_{j_1}+\cdots+y_{j_k} 
\quad\textup{for 
  all $1\leq k\leq n$ \quad and for all $1\leq j_1<\ldots<j_k\leq n$}
\end{equation}
We say that $\xul$ {\em majorizes\/} $\yul$ and write
$\xul\majgeq\yul$ if   $\xul\wmajgeq\yul$ and
$x_1+\ldots+x_n=y_1+\ldots+y_n$.

Observe the following:   
for real $n$-tuples $\xul$ and $\yul$ with
$\xul\majgeq\yul$,   given any $k$, $1\leq k\leq n$,  and any sequence $1\leq
i_1<\ldots<i_k\leq n$,    we have 
\begin{equation}\label{e:majgen}
y_{i_1}+\cdots+y_{i_k}\ \geq\ x_{n-k+1}^\downarrow+\cdots+x_n^\downarrow
\end{equation}
Indeed, let
$\{i_{k+1},\ldots,i_n\}:=\{1,\ldots,n\}\setminus\{i_1,\ldots,i_k\}$.
Then, on the one hand, when $y_{i_{k+1}}+\cdots+y_{i_n}$ is added to
the left hand side and 
$\xdown_1+\cdots+\xdown_{n-k}$ to the 
right hand side the resulting quantities are equal,  and,  on the 
other,   $y_{i_{k+1}}+\cdots+y_{i_n}\leq
\ydown_1+\cdots+\ydown_{n-k}\leq\xdown_1+\cdots\xdown_{n-k}$. 
\mysubsection{Majorization and weak interlacing}\mylabel{ss:majwil}
Let $\lseq:\lambda_1\geq\ldots\geq\lambda_n$ be a non-decreasing
sequence of real numbers and
$\mseq=(\mu_1,\ldots,\mu_n)\in\mathbb{R}^n$ such that
$\lseq\majgeq\mseq$.
Then $\lseq\winterlace(\mu_1,\ldots,\mu_{n-1})$.    (Proof:
(\ref{e:d:maj:cor}) and (\ref{e:majgen}).)

\mysubsection{Partitions}\mylabel{ss:partitions}  A {\em partition\/}
is a non-increasing sequence of non-negative integers that is
eventually zero.   Example:
$6$, $4$, $4$, $3$, $1$, $0$, $0$, \ldots.   The non-zero elements of the
sequence are called the {\em parts\/} of the partition.   If the sum
of the parts of a partition $\piseq:\pi_1\geq\pi_2\geq\ldots$ is $n$,  the partition is said to be a
{\em partition of~$n$\/}, and we write $|\piseq|=n$.   The example above is a partition of $18$
with $5$ parts.   

The trailing zeros in a partition are non-significant.  Thus $6$, $4$,
$3$, $3$, $1$, $0$, $0$, \ldots may also be written as $6$, $4$, $3$,
$3$, $1$.   

Each partition has an associated {\em shape\/}.
Given a partition
$\piseq: \pi_1\geq\pi_2\geq\ldots$ of $n$, its associated shape consists of a grid of $n$ squares, all of the same size, arranged top- and left-justified, with $\pi_1$~squares in the first row,  $\pi_2$ squares in the second, and so on (the rows are counted from the top downwards).     The shape corresponding to the partition $6$, $4$, $3$, $3$, $1$,  for example,  is this:
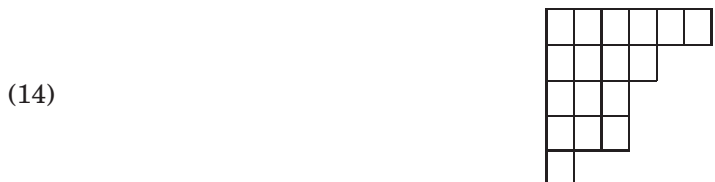
\begin{equation}\label{e:shape:5322}\begin{array}{|c|c|c|c|c|c|}
\hline 
& & &  & &\\ 
\hline 
& & &\\
\cline{1-4}
& &\\
\cline{1-3}
& &\\
\cline{1-3}
 \\
\cline{1-1}
\end{array}\end{equation}

We say that a partition {\em fits into a rectangle $(a,b)$\/},  where
$a$ and $b$ are non-negative integers,  if the number of parts is at
most $a$ and every part  is at most $b$.   The
terminology should make sense if we think of the shape associated to a
partition.   The partition whose shape is displayed above fits into
the rectangle $(a,b)$ if and only if $a\geq 5 $ and $b\geq 6$.
\mysubsubsection{Complementary partitions}\mylabel{sss:comppart}
Let $\piseq: \pi_1\geq\pi_2\geq\ldots$ be a partition that fits
into the rectangle $(a,b)$---in other words,  $b\geq
\pi_1$ and $\pi_j=0$ for $j>a$.    The {\em complement to~$\piseq$ in the rectangle
$(a,b)$\/} is the partition $\picseq$ defined as follows:
$\pic_j=b-\pi_{a+1-j}$ for $1\leq j\leq a$ and  $\pic_j=0$ for $j>a$.
For example,
the complement of the partition $6$, $4$, $3$, $3$, $1$  in
the rectangle $(7,6)$ is $6$, $6$, $5$, $3$, $3$, $2$.
\mysubsection{Colored partitions}\mylabel{ss:cpart}
Let $r$ be a positive integer.    An {\em $r$-colored partition\/} or
a {\em partition into $r$ colors\/} is a partition in which
each part is assigned an integer between $1$ and $r$.   The number
assigned to a part is its {\em color\/}.    We may think of an
$r$-colored partition as just an ordered $r$-tuple of
$(\piseq^1, \ldots, \piseq^r)$ of partitions:   the partition $\piseq^j$ consists
of all parts of color $j$ of the $r$-colored partition.    An
{\em $r$-colored partition of~$n$\/} is an $r$-colored partition with
$|\piseq^1|+\cdots+|\piseq^r|=n$.
\mysubsection{Partition overlaid patterns}\mylabel{ss:pop} A {\em
  partition overlaid pattern\/} ({\em POP\/} for short) consists of an
integral GT pattern $\lseq^1$, \ldots, $\lseq^n$,  and, for every
ordered pair $(j,i)$ of integers with $1\leq j<n$ and $1\leq i\leq j$, a partition $\pijiseq$ that fits into the rectangle $(\lambda^{j+1}_i-\lambda^j_i,\lambda_i^j-\lambda^{j+1}_{i+1})$.   Example: a partition overlay on the pattern displayed in (\ref{e:patex}) consists of three partitions $\pitwoseq^1$, $\pitwoseq^2$, and $\pioneseq^1$ that fit respectively into the rectangles $(0,2)$, $(1,1)$, and $(2,1)$. 

POPs 
parametrize bases of local Weyl modules of current algebras of type $A$
(as observed in~\S\ref{ss:clbase} below) just
as integral GT patterns parametrize bases of irreducible
representations of simple Lie algebras of type $A$ (as proved by Gelfand-Tsetlin and is well
known).

The {\em bounding sequence\/}, {\em area\/}, {\em
  trapezoidal area\/}, {\em weight\/}, etc.\ of a POP are just the
corresponding notions attached to the underlying pattern.
The {\em number of boxes\/} in a POP~$\pop$ is the sum
$\sum_{(j,i)}|\pijiseq|$  of the number of boxes in each of its
constituent partitions.    It is denoted by~$|\pop|$.   
Among POPs with a fixed underlying pattern,   the maximum possible
value of the number of boxes is evidently the area of the
pattern.    The {\em depth\/} of a POP~$\pop$  is defined by $\depth{\pop}:=
\traparea{\pattern}-|\pop|$,  where $\pattern$ is
the underlying pattern of~$\pop$.
\mysubsection{Weights identified as tuples}\mylabel{ss:wttup}
Let $\lieg=\slrpone$ be the simple Lie algebra consisting of $(r+1)\times(r+1)$ complex
traceless matrices ($r\geq1$).      
Let $\lieh$ and $\lieb$ be respectively the diagonal and upper
triangular subalgebras of~$\lieg$.    Linear functionals on~$\lieh$
are called {\em weights\/}.

Let $\epsilon_i$, $1\leq i\leq r+1$, be the weight that maps a
diagonal matrix to its entry in position~$(i,i)$.
Observe that $\epsilon_1+\cdots+\epsilon_{r+1}=0$.
Every weight may be expressed as
$a_1\epsilon_1+\cdots+a_{r+1}\epsilon_{r+1}$, with
$\aseq\in\crpone$.   Two elements 
in~$\crpone$ are said to be {\em equivalent\/} if their
difference is a multiple of $\one:=(1,\ldots,1)$, so that  weights are 
identified with equivalence classes in~$\crpone$.    

We will use this identification, often tacitly.  For a weight $\eta$, we denote by
$\underline{\eta}$ an element in the corresponding
equivalence class in~$\crpone$.    Depending upon the context,  this
$\etaseq$ may denote a particular representative:   we will see two
instances of this below.

A weight is {\em integral\/} if there exists a tuple~$\aseq$ in~$\crpone$
consisting of integers that corresponds to it;   it is {\em dominant\/} if $a_1\geq\ldots\geq
a_{r+1}$.   These notions correspond to the respective notions in the representation
theory of~$\lieg$ (with respect to the choice of $\lieh$ as Cartan
subalgebra and $\lieb$ as Borel subalgebra).  Dominant integral weights are thus in bijection with 
integer tuples of the form
$\lambda_1\geq\ldots\geq\lambda_r\geq\lambda_{r+1}=0$.
As an example,  consider the highest root $\theta$ of~$\lieg$.   The corresponding
element of $\crpone$ is $\tseq=(2,1,\ldots,1,0)$.

Let $\lambda$ be a dominant integral weight and $V(\lambda)$ the
irreducible representation of~$\lieg$ with highest weight~$\lambda$.
An integral weight $\mu$ is a weight of 
$V(\lambda)$ if
and only if $\lseq\majgeq\mseq$,  where the tuple
$\mseq=(\mu_1,\ldots,\mu_{r+1})$ representing $\mu$ is so chosen that
$\lambda_1+\cdots+\lambda_{r+1}=\mu_1+\cdots+\mu_{r+1}$.   (See also
Proposition~\ref{p:wtmajbseq} and Theorem~\ref{t:main}
in~\S\ref{s:theorem} below.)

Fix an invariant form $(\ |\ )$ on $\lieh^\star$ such that for every
root $\alpha$ we have $(\alpha|\alpha)=2$.  Given
$\lambda\in\lieh^\star$, how do we compute $(\lambda|\lambda)$ in
terms of the corresponding tuple $\lseq$?    We have
$(\lambda|\lambda)=||\lseq||^2:=\lambda_1^2+\cdots+\lambda_{r+1}^2$
provided that $\lseq$ is so chosen that
$\lambda_1+\cdots+\lambda_{r+1}=0$.  We will have occasion to compute
$(\lambda|\lambda)-(\mu|\mu)$ for $\lambda$, $\mu$ in $\lieh^\star$.
We observe that it equals $||\lseq||^2-||\mseq||^2$ provided that
$\lseq$ and $\mseq$ satisfy $\lambda_1+\cdots+\lambda_{r+1}=\mu_1+\cdots+\mu_{r+1}$.

\mysection{On area maximizing Gelfand-Tsetlin patterns}\mylabel{s:theorem}\mylabel{s:area}
\noindent
This section is elementary and combinatorial.     Its purpose is to
prove Theorem~\ref{t:main} below.    The representation theoretic
relevance of the theorem is discussed in \S\ref{s:relevance}.   For an $n$-tuple
$\underline{x}:=(x_1,\ldots,x_n)$ of real numbers,  the norm is defined as usual:  $||x||:=\sqrt{x_1^2+\cdots+x_n^2}$.

We begin with a proposition which should be well known.    We state and prove it in
order to put things in context and in the interest of completeness.
\begin{proposition}\mylabel{p:wtmajbseq}
Let $\lseq^1$, \ldots, $\lseq^n$ be a GT pattern with weight $\mseq$.
Then $\lseq^n\majgeq\mseq$.
\end{proposition}
\bmyproof
Proceed by induction on~$n$.    In case $n=1$,  we have
$\mseq=\lseq^1$, and the result is obvious.   

Now suppose that $n\geq2$.   
Let  $\mseq^{n-1}=(\mu_1, \ldots, \mu_{n-1})$ 
denote the weight of $\lseq^1$, \ldots, $\lseq^{n-1}$.    
By the induction hypothesis, we have
$\lseq^{n-1}\majgeq\mseq^{n-1}$.  
Since $\lseq^{n-1}={\lseq^{n-1}}^\downarrow$,   this means the
following:    for any $k$, $1\leq k\leq n-1$,  and any sequence 
$1\leq i_1<\ldots<i_k\leq n-1$,   we have:
\begin{equation}\label{e:wtmajbseq:ind}
\lambda^{n-1}_1+\lambda_2^{n-1}+\cdots+\lambda_{k}^{n-1} 
\quad\geq\quad \mu_{i_1}+\cdots+\mu_{i_k}
\end{equation}

Since
$\lseq^n={\lseq^n}^\downarrow$,   in order to show $\lseq^n\majgeq\mseq$,
we need to prove the following two sets of inequalities: 
for any $k$, $0\leq k\leq n-1$,  and any sequence 
$1\leq i_1<\ldots<i_k\leq n-1$:    
\begin{equation}\label{e:wtmajbseq:1}
\lambda^n_1+\lambda_2^n+\cdots+\lambda_{k}^n\quad\geq\quad \mu_{i_1}+\cdots+\mu_{i_k}
\end{equation}
\begin{equation}\label{e:wtmajbseq:2}
\lambda^n_1+\lambda_2^n+\cdots+\lambda_{k+1}^n\quad\geq\quad \mu_{i_1}+\cdots+\mu_{i_k}+\mu_n 
\end{equation}
Equation (\ref{e:wtmajbseq:1}) follows by combining (\ref{e:wtmajbseq:ind})
with the inequalities $\lambda_1^n\geq \lambda_1^{n-1}$, \ldots,
$\lseq^n_k\geq\lseq^{n-1}_k$  (which hold since $\lseq^{n-1}\interlace\lseq^n$).
As to (\ref{e:wtmajbseq:2}),   since 
$\mu_n=(\lambda^n_1+\cdots+\lambda^n_n)-(\lambda^{n-1}_1+\cdots+\lambda^{n-1}_{n-1})$,
it is equivalent to:
\begin{equation}\label{e:wtmajbseq:remainder:1}
\lambda^{n-1}_1+\lambda_2^{n-1}+\cdots+\lambda_{k}^{n-1} \ +\ 
(\lambda^{n-1}_{k+1}-\lambda^n_{k+2})+\cdots+(\lambda^{n-1}_{n-1}-\lambda^n_n)\quad\geq\quad \mu_{i_1}+\cdots+\mu_{i_k}
\end{equation}
But each of
$(\lambda^{n-1}_{k+1}-\lambda^n_{k+2})$, \ldots,
$(\lambda^{n-1}_{n-1}-\lambda^n_n)$ is non-negative (because
$\lseq^{n-1}\interlace\lseq^n$),   and thus (\ref{e:wtmajbseq:remainder:1})
too follows from~(\ref{e:wtmajbseq:ind}).
\emyproof

Here is the main result of this section:
\begin{theorem}\mylabel{t:main}
Let $\lseq=\lseqlong$ be a non-increasing 
sequence of real numbers and $\mseq=(\mu_1,\ldots,\mu_n)$ an element of $\mathbb{R}^n$ that is majorized by $\lseq$:
$\lseq\majgeq\mseq$.    
Then there is a unique GT pattern $\pattern:\lseq^1,\ldots,\lseq^n$ with bounding 
sequence $\lseq^n=\lseq$, weight $\mseq$, and satisfying the following:
\begin{equation}\begin{split}\label{e:theorem}
&\textup{For any 
$j$, $1\leq j\leq n$,  its $j^\textup{th}$ row  $\lseqj$  
majorizes the $j^\textup{th}$ row $\kseqj$ of 
any pattern with}\\[-2mm]
&\textup{bounding sequence $\lseq$ and weight $\mseq$: \quad 
$\lseqj\majgeq\kseqj$.}
\end{split}\end{equation}
This unique pattern $\pattern$ has the following properties:
\begin{enumerate}[(A)]
\item It is integral if $\lseq$ and $\mseq$ are integral.    
\item\label{item:b} Its area equals
  $\frac{1}{2}(||\lseq||^2-||\mseq||^2)$, which is strictly more than
  the area of any other pattern with bounding sequence~$\lseq$ and weight~$\mseq$.
\end{enumerate}
\end{theorem}
\noindent
Before turning to the proof,  we make a couple of remarks.     
Consider the GT polytope consisting of patterns with bounding sequence $\lseq$ and
weight $\mseq$. 
   Property~(\ref{item:b}) says that the pattern~$\pattern$
is the unique solution to the problem of maximizing the area function
on this polytope.
It is
clear from~(\ref{e:theorem}) that~$\pattern$ is a vertex of the
polytope.
   Observe that, in general, maxima of quadratic
functions on polytopes are neither unique nor always vertices.

For the proof of the theorem,  which appears in~\S\ref{sss:pftmain}
below,  we now make preparations.
\begin{proposition}\mylabel{p:traparea}  Let $n\geq1$ and $\lseq$, $\lpseq$ be sequences of lengths
  $n$ and $n-1$ respectively that are interlaced:
  $\lseq\interlace\lpseq$.         Then the trapezoidal area
  $\traparea{\lseq,\lpseq}$ is given by
\begin{equation}\label{e:traparea}
2  \traparea{\lseq,\lpseq} = ||\lseq||^2-||\lpseq||^2-((\lambda_1+\cdots+\lambda_n)-(\lambda'_1+\cdots+\lambda'_{n-1}))^2 
\end{equation}
\end{proposition}
\bmyproof
 Proceed by induction on~$n$.   In the case $n=1$,  both sides vanish
 (when correctly interpreted).    
%
Now suppose that $n\geq 2$.    Let $\kseq$ and $\kpseq$ be the
sequences of length $n-1$ and $n-2$ obtained by deleting respectively
$\lambda_n$ from $\lseq$ and $\lambda'_{n-1}$ from $\lpseq$:  the
sequence $\kpseq$ is empty in case $n=2$.     Also, let
\begin{equation}
\linearn:=(\lambda_n+\sum_{j=1}^{n-1}(\lambda_j-\lambda'_j))\quad\quad\quad\quad
\linearnmone:=(\lambda_{n-1}+\sum_{j=1}^{n-2}(\lambda_j-\lambda'_j))
\end{equation}
so that the last term in the desired equation (\ref{e:traparea}) is
$-\linearn^2$.

Then, firstly, by
induction:
\begin{equation}\label{e:traparea:k}
2  \traparea{\kseq,\kpseq} = ||\kseq||^2-||\kpseq||^2-\linearnmone^2
\end{equation}
Secondly, as is easily seen:
\begin{align}\label{e:traparea:lk:one}
\traparea{\lseq,\lpseq}&=\traparea{\kseq,\kpseq}+
(\lambda'_{n-1}-\lambda_n)\sum_{j=1}^{n-1} (\lambda_j-\lambda'_j)\\
\label{e:traparea:lk}
&=\traparea{\kseq,\kpseq}+
(\lambda'_{n-1}-\lambda_n) (\linearnmone-\lambda'_{n-1})
\end{align}
Finally, since $\linearn=\linearnmone-(\lambda_{n-1}'-\lambda_n)$:
\begin{equation}\label{e:linearn:linearnmone}
\linearn^2=\linearnmone^2+(\lambda_{n-1}'-\lambda_n)^2-2(\lambda_{n-1}'-\lambda_n)\linearnmone
\end{equation}
Adding twice of (\ref{e:traparea:lk}) with (\ref{e:traparea:k}) and
(\ref{e:linearn:linearnmone}),     we get
\begin{align*}
2\traparea{\lseq,\lpseq}+\linearn^2&=||\kseq||^2-||\kpseq||^2-(\lambda'_{n-1}-\lambda_n)(\lambda_{n-1}'+\lambda_n)\\
&=(||\kseq||^2+\lambda_n^2)-(||\kpseq||^2+{\lambda'_{n-1}}^2)\\
&=||\lseq||^2-||\lpseq||^2 
\end{align*}
and the proposition is proved. \emyproof
\bcor\mylabel{c:traparea}
The trapezoidal area of a pattern $\pattern$: $\lseq^1$, $\lseq^2$, \ldots,
$\lseq^n$ is given by
\begin{equation}\label{e:traparea:pat}
\traparea{\pattern}=\frac{1}{2}(||\lseq||^2-||\mseq||^2) 
\end{equation}
where $\lseq=\lseq^n$ is the bounding sequence and $\mseq$ the weight
of $\pattern$.
\ecor
\bmyproof    By the proposition:
\begin{equation}\label{e:traprea:cor:proof}
\traparea{\lseq^j,\lseq^{j-1}}=\frac{1}{2}(||\lseq^j||^2-||\lseq^{j-1}||^2-\mu_j^2)
\quad\quad\textup{for $j=n, n-1, \ldots, 1$}
\end{equation}
Adding these $n$ equations gives us the desired result.
\emyproof
\blemma\mylabel{l:forrtpf}
Let $\lseq$ and $\lpseq$ be as in Proposition~\ref{p:traparea}.
Suppose further that
\begin{enumerate}
\item\label{i:one} 
$\lseq$ is integral;
\item\label{i:two}  $\lambda'_1+\cdots+\lambda_{n-1}'$ is an
integer; and
\item\label{i:three} $\traparea{\lseq,\lpseq}=\area{\lseq,\lpseq}$.
\end{enumerate}
Then $\lpseq$ is integral.
\elemma
\bmyproof 
Let $k$ be the largest integer, $1\leq k\leq n-1$, such that 
$\lambda'_j=\lambda_j$
for all $j<k$.     
From item~(\ref{i:three}) it follows that  $\lambda'_j=\lambda_{j+1}$ for all $j>k$. 
From item~(\ref{i:one}) it follows that $\lambda'_j$ is an integer for $j\neq k$.   
From item~(\ref{i:two}) it follows that $\lambda'_k$ is also an integer.
\emyproof
\bcor\mylabel{c:l:forrtpf}
Let $\lseq:\lambda_1\geq\ldots\geq\lambda_n$ be a non-increasing sequence of integers.  Let $\mseq$ be in~$\mathbb{Z}^n$ such that $\lseq\majgeq\mseq$. Let $\pattern$ be a GT pattern with bounding sequence $\lseq$, weight~$\mseq$, and
area $\frac{1}{2}(||\lseq||^2-||\mseq||^2)$.  Then $\pattern$ is integral.
\ecor
\bmyproof
Let $\lseqj$ denote the $j^{\textup{th}}$ row of $\pattern$.   By Corollary~\ref{c:traparea},   $\traparea{\pattern}=\area{\pattern}$,   so $\traparea{\lseqj,\lseqjmone}
=\area{\lseqj,\lseqjmone}$ for all $j$, $n\geq j\geq 1$.  Since $\mu_j=(\lambdaj_1+\cdots+\lambdaj_j)-(\ljmone_1+\cdots+\ljmone_{j-1})$,  it follows (by an easy decreasing induction) that $\lambdaj_1+\cdots+\lambdaj_j$ is an integer for all $j$, $n\geq j\geq 1$.  By applying Lemma~\ref{l:forrtpf} repeatedly, we see successively that $\lseq^{n-1}$, \ldots, $\lseq^1$ are all integral.
\emyproof

\mysubsection{Proof of 
  Theorem~\ref{t:main}}\mylabel{ss:thproof}
\blemma\mylabel{l:main}
Let $\lseq^n:\lambda^n_1\geq\ldots\geq\lambda^n_n$ be a non-increasing
sequence 
of real numbers and $\mseq=(\mu_1,\ldots,\mu_n)$ an element of~$\mathbb{R}^n$ such
that $\lseq^n\majgeq\mseq$.   Then there exists a unique non-increasing sequence
$\lseqnmone: \lambdanmone_1\geq \ldots\geq\lambdanmone_{n-1}$ of 
real numbers such that the following hold:
\begin{enumerate}
\item\label{i:one:l:main} $\lseqnmone\ \interlace\ \lseqn$
\item $\lambdanmone_1+\cdots+\lambdanmone_{n-1}\ =\ \mu_1+\cdots+\mu_{n-1}$
\item Let $\kseqn: \kappan_1\geq\ldots\geq\kappan_n$ be a 
  non-increasing sequence of length $n$ of real numbers and 
  $\kseqnmone=(\knmone_1,\ldots,\knmone_{n-1})$ an element 
  of~$\mathbb{R}^{n-1}$ such that:
\begin{quote}
(i)~$\lseqn\majgeq\kseqn$, \quad 
(ii)~$\kseqn\winterlace\kseqnmone$, \quad and \quad 
(iii)~$\knmone_1+\cdots+\knmone_{n-1}=\mu_1+\cdots+\mu_{n-1}$. 
\end{quote}
Then $\lseqnmone\majgeq\kseqnmone$. 
\end{enumerate}

\noindent
Moreover, the unique sequence $\lseqnmone$ has the following properties:
\begin{enumerate}[(a)]
\item $\lseqnmone$ is integral if $\lseqn$ and $\mseq$ are integral.
\item  $\traparea{\lseqn,\lseqnmone}=\area{\lseqn,\lseqnmone}$.
\item Let $\kseqnmone:
   \kseqnmone_1\geq\ldots\geq\kseqnmone_{n-1}$ be a non-increasing
   sequence of real numbers such that:
\begin{quote} (i')~$\kseqnmone\interlace\lseqn$,
   (ii')~$\knmone_1+\cdots+\knmone_{n-1}=\mu_1+\cdots+\mu_{n-1}$,  and
   (iii')~$\lseqnmone\neq\kseqnmone$.\end{quote}
Then  
$\traparea{\lseqn,\kseqnmone}\gneq\area{\lseqn,\kseqnmone}$.
\end{enumerate}
\elemma
\bmyproof     The uniqueness of $\lseqnmone$ is easy to see.    Indeed,
if $\etanmone$ be a another sequence with properties (1)--(3),   then
by applying (3) with $\kseqn=\lseq^n$ and $\kseqnmone=\etanmone$,  we
see that $\lseqnmone\majgeq\etanmone$.    By the same argument with
the roles of $\etanmone$ and $\lseqnmone$ switched,   we see that
$\etanmone\majgeq\lseqnmone$.        It follows from the definition
of~$\majgeq$ that it is a partial order, so
we conclude $\etanmone=\lseqnmone$.

We now turn to the existence of~$\lseqnmone$.
Consider the auxiliary non-decreasing sequence of terms  $\ltseq:
\lambdatn_1\leq\ldots\leq\lambdatn_n$,   where 
\begin{equation}\label{e:ltseqdef}
\lambdatn_k\ :=\  (\lambdan_1+\cdots+\lambdan_n)\ -\ \lambdan_k
\end{equation}
Since ${\lseqn}^\downarrow=\lseqn$,  it follows from equations
(\ref{e:majgen}) and (\ref{e:d:maj:cor})     that
\begin{equation}\label{e:lemma:mu:1}
  \lambdatn_1\ \leq\ \mu_1+\cdots+\mu_{n-1}\ \leq\ \lambdatn_n
\end{equation}
Fix a $\knot$, $1\leq \knot\leq n-1$, such that  
\begin{equation}\label{e:lemma:mu:2}
\lambdatn_\knot\ \leq\ 
    \mu_1+\cdots+\mu_{n-1}\ \leq\ \lambdatn_{\knot+1}\end{equation}
In fact, there is a unique such $\knot$ except when
$\mu_1+\cdots+\mu_{n-1}=\lambdatn_j$ for some $j$.    

Set
\begin{equation}\label{e:lmonej}
\lambdanmone_j:=\left\{ 
\begin{array}{ll}
\lambdan_j & \textup{for $j<\knot$}\\
\lambdan_{j+1} & \textup{for $j>\knot$}\\
\end{array}
\right. 
\end{equation}
\begin{equation}\label{e:lambdanmonek}
\lambdanmone_{\knot}:=
(\mu_1+\cdots+\mu_{n-1})-(\lambdan_1+\cdots+\lambdan_{\knot-1}+\lambdan_{\knot+2}+\cdots+\lambdan_n)
\end{equation}
From (\ref{e:lambdanmonek}) and (\ref{e:lemma:mu:2}),   we see that 
\begin{equation}\label{e:lnillnmone}
\lambdan_\knot-\lambdanmone_{\knot} = \lambdatn_{\knot+1}-(\mu_1+\cdots+\mu_{n-1})
\geq 0
\quad\text{\&}\quad
\lambdanmone_\knot-\lambdan_{\knot+1} =
(\mu_1+\cdots+\mu_{n-1})-\lambdatn_\knot\geq 0 
\end{equation}
To observe that $\lseqnmone$ is a non-increasing sequence,   
first note that, from (\ref{e:lmonej}) and the fact that $\lseqn$ is non-increasing, we have:
\begin{equation}
\lambdanmone_1\geq
\ldots\geq\lambdanmone_{\knot-1}\quad\quad\quad
\lambdanmone_{\knot+1}\geq\ldots\geq\lambdanmone_{n-1}
\end{equation}
And, combining (\ref{e:lmonej}), (\ref{e:lnillnmone}), and the
non-increasing property of $\lseqn$, we have:
\begin{equation}
\lambdanmone_{\knot-1}=\lambdan_{\knot-1}\geq\lambdan_\knot\geq\lambdanmone_\knot
\quad\quad\quad
\lambdanmone_\knot\geq\lambdan_{\knot+1}\geq\lambdan_{\knot+2}=\lambdanmone_{\knot+1}
\end{equation}

Item~(1) follows from (\ref{e:lmonej}) and
(\ref{e:lnillnmone}),  item (2) from
(\ref{e:lmonej}) and (\ref{e:lambdanmonek}).    
As to item (3),   we must show that 
for any $j$, $1\leq j\leq 
n-1$,  and any choice $1\leq i_1<\ldots<i_j\leq n-1$, the following holds:
\begin{equation}\label{e:lnmone:a}
\lnmone_1+\cdots+\lnmone_j\quad\geq \quad
  \knmone_{i_1}+\cdots+\knmone_{i_j}
\end{equation}
For $j<\knot$,   we have,  by successively using (\ref{e:lmonej}),
(i), and (ii):
\begin{equation}
\lnmone_1+\cdots+\lnmone_j \quad=\quad
\lambdan_1+\cdots+\lambdan_j\quad\geq\quad
\kappan_1+\cdots+\kappan_j\quad\geq\quad
  \knmone_{i_1}+\cdots+\knmone_{i_j}
\end{equation}
For $j\geq \knot$,   substituting for 
$\lambdanmone_i$ from (\ref{e:lmonej}) and (\ref{e:lambdanmonek}),
the left hand side of (\ref{e:lnmone:a}) becomes:
\begin{align*}
\lnmone_1+\cdots+\lnmone_j\quad 
=\quad(\mu_1+\cdots+\mu_{n-1})-(\lambdan_{j+2}+\cdots+\lambdan_n) 
\end{align*}
Using (iii), we get
\begin{equation*}\lnmone_1+\cdots+\lnmone_j\quad =\quad
(\knmone_1+\cdots+\knmone_{n-1})-(\lambdan_{j+2}+\cdots+\lambdan_n) 
\end{equation*}
Letting $\{i_{j+1},\ldots,i_{n-1}\}$ denote the complement 
$\{1,\ldots,n-1\}\setminus\{i_1,\ldots,i_j\}$, the right hand side of
the last equation may be rewritten as 
\begin{equation*}
(\knmone_{i_1}+\cdots+\knmone_{i_j})\quad+\quad\left((\knmone_{i_{j+1}}+\cdots+\knmone_{i_{n-1}})-(\lambdan_{j+2}+\cdots+\lambdan_n) \right) 
\end{equation*}
The second parenthetical term here is non-negative,  for  by~(ii)
and~(i)
\begin{equation*}
\knmone_{i_{j+1}}+\cdots+\knmone_{i_{n-1}}\quad\geq\quad
\kappan_{j+2}+\cdots+\kappan_n\quad\geq\quad
\lambdan_{j+2}+\cdots+\lambdan_n
\end{equation*}
and the proof of item~(3) is complete.

Assertion (a) is immediate from the definitions (\ref{e:lmonej}) and
(\ref{e:lambdanmonek}) of~$\lseqnmone$.
As to (b),  it is clear from definitions (\ref{e:defarea}) 
and (\ref{e:deftraparea}) that it holds, 
since either $\lambdan_j=\lambdanmone_j$ or 
$\lambdanmone_j=\lambdan_{j+1}$ for every $j\neq\knot$. 
%
%
Towards the proof of (c), let $\ell$ and $\rho$ be respectively the smallest and largest $j$,
$1\leq j\leq n-1$,  such that $\lnmone_j\neq\knmone_j$.    Taking
$\kseqn$ to be $\lseqn$ in (3),  we obtain 
$\lseqnmone\majgeq\kseqnmone$, which implies that
$\lnmone_\ell>\knmone_\ell$ and $\knmone_\rho>\lnmone_\rho$.
In particular,  $\ell<\rho$.
From (\ref{i:one:l:main}),  we have
$\lambdan_\ell\geq\lnmone_\ell>\knmone_\ell$ and
$\knmone_\rho>\lnmone_\rho\geq\lambdan_{\rho+1}$,  so that 
$(\lambdan_\ell-\knmone_\ell)(\knmone_\rho-\lambdan_{\rho+1})$ is a
non-trivial contribution to $\traparea{\lseqn,\kseqnmone}-\area{\lseqn,\kseqnmone}$.
\emyproof
\bcor\mylabel{c:lemma:main}
With hypothesis and notation as in the lemma,  $\lseqnmone\majgeq(\mu_1,\ldots,\mu_{n-1})$.
\ecor
\bmyproof
Put $\kseqn=\lseqn$ and $\kseqnmone=(\mu_1,\ldots,\mu_{n-1})$ in
item~(3).  Hypothesis~(ii) holds by the observation in~\S\ref{ss:majwil}.
\emyproof
\mysubsubsection{Proof of Theorem~\ref{t:main}}\mylabel{sss:pftmain}
The uniqueness of the pattern~$\pattern$ being obvious,  it is enough to
prove its existence.
Apply Lemma~\ref{l:main} to the given pair $\lseq$ and $\mseq$ (by
taking $\lseq^n$ in the statement of the lemma to be $\lseq$).    The
$\lseqnmone$ we obtain as a result is such that
$\lseqnmone\majgeq(\mu_1,\ldots,\mu_{n-1})$ 
(Corollary~\ref{c:lemma:main}) so we can apply the lemma
again,  this time to the pair $\lseqnmone$ and
$(\mu_1,\ldots,\mu_{n-1})$.   Continuing thus,   we obtain, by
items (1) and (2) of the lemma, a GT
pattern---let us denote it~$\pattern$---with bounding sequence
$\lseq$ and weight $\mseq$.

We claim that the pattern $\pattern$ 
satisfies~(\ref{e:theorem}).  To prove this, proceed by reverse induction on~$j$.    For $j=n$,
we have $\lseq^n=\kseq^n=\lseq$,  so the statement is evident. 
For the induction step, suppose we have proved that 
$\lseq^j\majgeq\kseq^j$.     Note that $\lseq^{j-1}$ is constructed 
by applying Lemma~\ref{l:main} with $\lseq^j$ in place of $\lseq^n$
(in the notation of the lemma).    Assertion~(\ref{e:theorem}) follows from 
item~(3) of the lemma,   by substituting respectively $\kseq^j$,
$\kseq^{j-1}$, and $(\mu_1,\ldots,\mu_{j-1})$ for $\kseq^n$,
$\kseq^{n-1}$, and $(\mu_1,\ldots,\mu_{n-1})$. 

The pattern $\pattern$ is integral if $\lseq$ and $\mseq$ are so
(Lemma~\ref{l:main}~(a)),  so~(A) is clear.

Finally we prove~(B). 
Let $\pattern'$ be a pattern distinct from~$\pattern$ with bounding 
sequence $\lseq$ and weight~$\mseq$.   By Corollary~\ref{c:traparea}, 
$\traparea{\pattern}=\traparea{\pattern'}=(||\lseq||^2-||\mseq||^2)/2$.
By Lemma~\ref{l:main}~(b),  $\traparea{\pattern}=\area{\pattern}$.
Let $j$ be largest, $1\leq j\leq 
n$, such that the $j^\textup{th}$ row~$\kseq^j$ of $\pattern'$ is 
distinct from the $j^\textup{th}$ row of~$\pattern$.    
Then $j<n$ and, by Lemma~\ref{l:main}~(c),
$\traparea{\kseqjpone,\kseqj}\gneq\area{\kseqjpone,\kseqj}$,
so $\traparea{\pattern'}\gneq\area{\pattern'}$, and (B) is proved. \hfill $\Box$

\mysection{Relevance of the main theorem to the theory of Local Weyl modules}\mylabel{s:relevance}
\noindent
In this section we discuss the relevance of our main theorem (Theorem~\ref{t:main}) of~\S\ref{s:theorem} to representation theory.  We do this by means of giving a representation theoretic proof of a version of the theorem: see~\S\ref{ss:rtproof} below.     The proof is based on the theory of local Weyl modules for current algebras.   We first recall the required results from this theory.   

The reference for most of the basic material recalled here is~\cite{kac}.   While our notation does not follow that of~\cite{kac} faithfully,   there should be no cause for confusion.

\mysubsection{The current algebra~$\currlieg$ and the affine algebra~$\ghat$}\mylabel{ss:currlieg}
Let $\lieg$ be a complex simple finite dimensional Lie algebra.   
The corresponding {\em current algebra\/}, denoted $\curralg$,   is
merely the extension of scalars to the polynomial ring $\mathbb{C}[t]$
of~$\lieg$.   In other words,  $\curralg=\lieg\tensor\mathbb{C}[t]$,
with Lie bracket $[X\tensor f,Y\tensor g]=[X,Y]\tensor fg$ for all
$X$, $Y$ in $\lieg$ and $f$, $g$ in $\mathbb{C}[t]$.
   There is a natural grading on $\curralg$ given by the degree in~$t$:   thus $X\tensor t^s$ has degree $s$ (here $X$ is in~$\lieg$ and $s$ is a non-negative integer).   There is an induced grading by non-negative integers on the universal enveloping algebra~$U(\currlieg)$.    We can talk about graded modules (graded by integers) of~$\currlieg$.      

Let $\ghat=\lieg\tensor\mathbb{C}[t,t^{-1}]\oplus\mathbb{C}c\oplus\mathbb{C}d$ be the affine algebra corresponding to~$\lieg$.  
The current algebra~$\currlieg$ is evidently a subalgebra of~$\ghat$.
\mysubsection{Fixing notation}\mylabel{ss:handb}
Fix a Cartan subalgebra~$\lieh$ of~$\lieg$ and a Borel subalgebra $\lieb$ containing $\lieh$ of~$\lieg$.    
Let $\lieg=\nminus\oplus\lieh\oplus\nplus$ be the triangular decomposition of $\lieg$ with $\lieb=\lieh\oplus\nplus$.   
Let $(\ |\ )$ denote the invariant form on $\lieh^\star$ such that $(\alpha|\alpha)=2$ for all long roots $\alpha$.

Put $\hhat:=\lieh\oplus\mathbb{C}c\oplus\mathbb{C}d$  and $\bhat:=\lieg\tensor t\mathbb{C}[t]\oplus\lieb\oplus\mathbb{C}c\oplus\mathbb{C}d$.     Denote by $\Lambda_0$ and $\delta$ elements of~$\hhat^\star$ 
such that $\langle \Lambda_0,c\rangle=\langle\delta,d\rangle=1$ and $\langle\Lambda_0,\lieh\rangle=
\langle\Lambda_0,d\rangle=\langle\delta,\lieh\rangle=\langle\delta,c\rangle=0$.    
Extend $(\ |\ )$ to a symmetric bilinear form on~$\hhat^\star$ by setting
$(\lieh^\star|\Lambda_0)=(\lieh^\star|\delta)=(\Lambda_0|\Lambda_0)=(\delta|\delta)=0$,   
$(\Lambda_0|\delta)=1$,  where $\lieh^\star$ is identified as the subspace of~$\hhat^\star$ that kills $c$ and $d$.

Fix a dominant integral weight $\lambda$ of~$\lieg$ (with respect to $\lieh$ and $\lieb$).     
\mysubsection{The local Weyl module~$\wml$}\mylabel{ss:lwm}  
An element $\hwel$ of a $\currlieg$-module is of {\em highest weight\/}~$\lambda$ if:
\begin{equation}
\label{e:hwt}
  \quad (\nplus\tensor\mathbb{C}[t])\,\hwel=0,\quad\quad (\lieh\tensor t\mathbb{C}[t])\,\hwel=0,\quad\quad H\,\hwel=\langle\lambda,H\rangle\hwel\quad\textup{ for $H\in\lieh$}  \end{equation} 
%
The {\em local Weyl module\/}~$\wml$ corresponding to $\lambda$ is the cyclic $\curralg$-module generated by an element $\hwel$ of highest weight~$\lambda$ (in other words, subject to the relations (\ref{e:hwt})) and further satisfying:
\begin{equation}
\label{e:wmtwo}  \lieg_{-\alpha}^{\langle\lambda,\alpha^\vee\rangle+1} \hwel=0 \textup{ for every simple root $\alpha$ of~$\lieg$}
\end{equation}
where $\alpha^\vee$ is the co-root corresponding to $\alpha$ and
$\lieg_{-\alpha}$ is the root space in $\lieg$ corresponding to $-\alpha$.
It is evident that $\wml$ is graded (since the relations in~(\ref{e:hwt}) and (\ref{e:wmtwo}) are all homogeneous).   We let the generator~$\hwel$ have grade~$0$,  so that $\wml=U(\nminus\tensor\mathbb{C}[t])\, \hwel$ is graded by the non-negative integers.
It is well known---the proofs are analogous to those in \cite[\S2]{cprt2001}---that $\wml$ is finite dimensional and moreover that it is maximal among finite dimensional modules generated by an element of highest weight~$\lambda$ (which means, more precisely, that for any finite dimensional cyclic $\currlieg$-module~$M$ generated by an element $m$ of highest weight~$\lambda$,  there exists a unique $\currlieg$-module map from $\wml$ onto $M$ mapping~$u$ to~$m$).
\mysubsection{Local Weyl modules as Demazure modules}\mylabel{ss:lwmdem}
Let $w_0$ be the longest element of the Weyl group~$W$ of $\lieg$. 
  Let $\Lambda$ be the dominant integral weight of~$\ghat$ and $w$ an element of the affine Weyl group such that
\begin{equation}\label{e:Lambda}
w\Lambda= t_{w_0\lambda}\Lambda_0  \quad\quad\textup{($w$ in the affine Weyl group, $\Lambda$ dominant)}
\end{equation}
where (see \cite[(6.5.2)]{kac}):
\begin{equation}\label{e:tgamma} 
t_\gamma\zeta:=\zeta+\langle \zeta,c \rangle\gamma-((\zeta|\gamma)+\frac{1}{2}(\gamma|\gamma)\langle 
\zeta, c\rangle)\delta  \quad\quad \textup{for $\gamma$ in $\lieh^\star$ and $\zeta$ in $\hhat^\star$}
\end{equation}
From the above two equations, we obtain that $\Lambda$ has level~$1$,  or, in other words $\langle \Lambda,c\rangle=1$:
\begin{equation}\label{e:level}
\langle \Lambda, c\rangle = (\Lambda|\nu(c))=(\Lambda|\delta)=(w\Lambda|w\delta)=(w\Lambda|\delta)=(t_{w_0\lambda}\Lambda_0|
\delta)=(\Lambda_0|\delta)=1
\end{equation}
where $\nu:\hhat\to\hhat^*$ is the isomorphism as
in \cite[\S2.1]{kac}:    $(w\Lambda|w\delta)=(w\Lambda|\delta)$
because of the invariance of the form under the action of the affine Weyl
group and the fact that $\delta$ is fixed by the affine Weyl group;
the penultimate equality holds because of (\ref{e:tgamma}). 

Let $V(\Lambda)$ be the integrable $\ghat$-module with highest weight~$\Lambda$.     With $w$ as above, denote by $\dem$ the Demazure submodule $U(\bhat)(V(\Lambda)_{w\Lambda})$ of $V(\Lambda)$:  here $V(\Lambda)_{w\Lambda}$ denotes the (one-dimensional) $\hhat$-weight space of $V(\Lambda)$ of weight $w\Lambda$.

We recall the identification of local Weyl modules as Demazure modules:
\begin{theorem}\mylabel{t:lwmdem} \textup{(\cite[1.5.1~Corollary]{cladv2006} for type~$A$ and \cite[Theorem 7]{fladv2007} in general)} Let $\lieg$ be simply laced and $\lambda$ be a dominant integral weight of~$\lieg$.  
With $w$ and $\Lambda$ as in~(\ref{e:Lambda}),   let $v$ be a non-zero element of the one-dimensional $\hhat$-weight space 
of the Demazure module $\dem$ of weight~$w_0w\Lambda= t_\lambda
w_0\Lambda_0=t_\lambda \Lambda_0$.    Then $v$ satisfies the relations (\ref{e:hwt}) and therefore, since $\dem$ is finite dimensional,  there exists a unique $\currlieg$-map from the local Weyl module~$\wml$ onto $\dem$ mapping the generator to $v$.    
This map is an isomorphism.  
\end{theorem}
\mysubsection{The key proposition}\mylabel{ss:keyprop}
We now state and prove the main proposition of this section.  
Observe that the $\lieg$-weights of the local Weyl module $W(\lambda)$ are
precisely the weights of the irreducible $\lieg$-module $V(\lambda)$
with highest weight $\lambda$:     indeed this follows from the finite
dimensionality of $\wml$ and the fact that $\wml=U(\nminus[t])u$,
where $u$ is the generator.
\begin{proposition}\mylabel{p:key}
Let $\lambda$ be a dominant integral weight of $\lieg$
and $\mu$ a weight of the irreducible $\lieg$-module $V(\lambda)$.   Let $M$ be the largest integer such that the $M^{\textup{th}}$-graded piece~$\wmlmM$ of the $\mu$-weight space of the local Weyl module $W(\lambda)$ is non-zero.  Then,  under the assumption that $\lieg$ is simple of simply laced type, we have:
\begin{enumerate}  
\item $M=\frac{1}{2}\cdot\left((\lambda|\lambda)-(\mu|\mu)\right)$
\item $\wmlmM$ has dimension~$1$
\end{enumerate}
\end{proposition}
\bmyproof
Fix notation as in~\S\ref{ss:lwmdem}.  Identify $\wml$ with the Demazure module $\dem$ as explained there.      
Since $\wml=U(\nminus[t])\cdot u$ and the generator $u$ is identified with an element of $\dem_{t_\lambda\Lambda_0}$,  the
$\hhat$-weights of $\wml$ are of the form \begin{equation}\label{e:wform}t_\lambda\Lambda_0-\eta+d\delta\end{equation}  where $\eta$ is a non-negative integral linear combination of the simple roots of~$\lieg$,  and $d$ is a non-negative integer.     
 From (\ref{e:tgamma}) we have:
\beq\label{e:tl0}
t_\lambda\Lambda_0=\Lambda_0+\lambda-\frac{1}{2}(\lambda|\lambda)\delta
\eeq
so we may rewrite (\ref{e:wform}) as:
\beq\label{e:wfrom2}
\lambda-\eta+\Lambda_0-\frac{1}{2}(\lambda|\lambda)\delta+d\delta
\eeq
Observe that this weight acts on $\lieh$ as $\lambda-\eta$.    Thus
$\wml_\mu$  is exactly the direct sum of $\hhat$-weight spaces of
$\dem$ corresponding to weights of the form (\ref{e:wform}) with $\mu=\lambda-\eta$.

Let $\eta$ be such that  $\lambda-\eta=\mu$.    
 Then by the hypothesis of maximality of~$M$, we have:
 \begin{itemize}
\item $\kappa:=t_\lambda\Lambda_0+(\mu-\lambda)+M\delta$ is a weight
 of $\dem$ 
\item but $\kappa_1:=t_\lambda\Lambda_0+(\mu-\lambda)+M'\delta$ is not
 (a weight of $\dem$) for $M'>M$.\end{itemize}   Then~$\kappa_1$
 is not a weight of $V(\Lambda)$ either  (this should be well known;
for a proof in our specific case, see~\cite[Lemma in~\S3.3.5]{br:thesis}).
Thus $\kappa$ is a {\em maximal weight\/} of~$V(\Lambda)$ in the sense of~\cite[\S12.6]{kac}.    Since $\Lambda$ is of level one (\ref{e:level}),  there exists, by~\cite[Lemma~12.6]{kac},   an element~$\gamma$ of the root lattice of~$\lieg$  such that $\kappa=t_\gamma \Lambda$.
In particular, $\kappa$ is an affine Weyl group translate of the highest weight~$\Lambda$ of~$V(\Lambda)$,    and so the multiplicity of the $\kappa$-weight space in the Demazure module $\dem$ (and so also in~$\wml$) cannot exceed~$1$.  This proves (2).

We now prove (1), by equating two expressions for $\kappa$.   On the one hand,  by the definition of $\kappa$ and (\ref{e:tl0}), we get
\beq\label{e:k1}
\kappa=\Lambda_0+\mu-\frac{1}{2}(\lambda|\lambda)\delta+M\delta
\eeq   
On the other hand we have $\kappa=t_\gamma\Lambda$.
Since $\Lambda$ is of level~$1$,  we obtain using \cite[(6.5.3)]{kac} that
\newcommand\barLam{\overline{\Lambda}}
\beq\label{e:tgL}
\kappa= t_\gamma\Lambda=\Lambda_0+(\barLam+\gamma)+\frac{1}{2}((\Lambda|\Lambda)-(\barLam+\gamma|\barLam+\gamma))\delta
\eeq
where $\barLam$ denotes the restriction to $\lieh$ of $\Lambda$. 
We have $(\Lambda|\Lambda)=(w\Lambda|w\Lambda)=(t_{w_0\lambda}\Lambda_0|t_{w_0\lambda}\Lambda_0)=0$,  where the last equality follows from (\ref{e:tl0}).

Equating the $\lieh^\star$ components on the right hand sides of (\ref{e:k1}) and (\ref{e:tgL}),  we get $\mu=\barLam+\gamma$;    now equating the coefficients of~$\delta$,   we get~(1).
\emyproof

Combining the proposition above with a result of Kodera-Naoi~\cite[\S3]{kn},  one can recover the known 
fact that the space of current algebra homomorphisms between two local Weyl modules is at most one dimensional: see \cite[Corollary~3.3.5]{br:thesis}.
\mysubsection{Bases for local Weyl modules in
type~$A$}\mylabel{ss:clbase} In this section, we construct certain
bases for local Weyl modules of the current algebra $\slrpone[t]$. These bases are indexed by POPs (see \S\ref{ss:pop}), and turn out to coincide, up to scaling, with the bases constructed by Chari and Loktev in \cite{cladv2006}.

\mysubsubsection{}\mylabel{sss:cl}    Let $\lieg=\slrpone$. Fix notation and terminology as in~\S\ref{ss:wttup}.
Fix a dominant integral weight $\lambda$ with corresponding tuple $\lseq:\lambda_1\geq\ldots\geq\lambda_r\geq\lambda_{r+1}=0$.
Let $x_{ij}^-$,  for $1\leq i\leq j\leq r$, denote the $(r+1)\times(r+1)$ complex matrix all of whose entries are
zero except the one in position~$(j+1,i)$ which is~$1$.   

Let $a$, $b$ be non-negative integers and $\piseq$ a partition that fits into
the rectangle $(a,b)$.   For~$k$ an integer with $1\leq k\leq b$,  let
$\mynum(a,b,\piseq,k)$ denote the number of parts of $\piseq$ that
equal~$k$.    Let $\mynum(a,b,\piseq,0)$ denote
$a-\sum_{k=1}^b\mynum(a,b,\piseq,k)$.   Let $x_{ij}^-(a,b,\piseq)$
denote
\begin{equation}\label{e:xijpop}
(x_{ij}^-\tensor1)^{(\mynum(a,b,\piseq,0))}\cdot
(x_{ij}^-\tensor t^1)^{(\mynum(a,b,\piseq,1))}\cdot\ \cdots\ \cdot
(x_{ij}^-\tensor t^b)^{(\mynum(a,b,\piseq,b))}
\end{equation}
where, for an operator $T$ and a non-negative integer $n$,  the symbol
$T^{(n)}$ denotes the {\em divided power\/} $T^n/n!$.
The order of factors in the above product is immaterial since they
commute with each other, so we may simply write $x_{ij}^-(a,b,\piseq)= \prod_{k=0}^b(x_{ij}^-\tensor t^k)^{(\mynum(a,b,\piseq,k))}$.

Let~$\pop$ be a POP with bounding sequence $\lseq$.  Let
$\loneseq$, \ldots, $\lrseq$, $\lrponeseq=\lseq$ be the rows of
$\pop$,  and $\pijiseq$, $1\leq i\leq j\leq r$, be the partition
overlay.    For $1\leq j\leq r$,  denote by $\rho^j_\pop$ the monomial
\[\prod_{i=1}^j
x_{ij}^-(\lambda_i^{j+1}-\lambda^j_i,\lambda_i^j-\lambda^{j+1}_{i+1}, \pijiseq)\]
where again the order is immaterial in the product.   Now define the
following:   
\begin{equation}\label{e:clmonom}
\rhocl{\pop}:=\rhocl{\pop}^1\cdot\rhocl{\pop}^2\cdot\ \cdots\ \cdot\rhocl{\pop}^{r-1}\cdot\rhocl{\pop}^r \, 
\quad\quad\quad
\vcl{\pop}:=\rhocl{\pop} w_\lambda\end{equation}
The order of the factors matters  in the expression for $\rhocl{\pop}$.

We shall call $\rhocl{\pop}$ the
 {\em Chari-Loktev (or CL)} monomial corresponding to~$\pop$ and, in
 view of the following theorem, $\vcl{\pop}$ the
 {\em Chari-Loktev (or CL)} basis element corresponding to~$\pop$.
\bthm\mylabel{t:cl}    The elements
$v_\pop$,  as $\pop$ varies over all POPs with
bounding sequence~$\lseq$,  form a basis for the local Weyl module~$W(\lambda)$.
\ethm
\bmyproof
We will show that these basis elements are scalar multiples of the
basis elements constructed by Chari-Loktev in \cite{cladv2006}. To
this end, we first recall the Chari-Loktev construction from \cite[\S
2]{cladv2006}. Let $m_i = \lambda_i-\lambda_{i+1}$ for $1 \leq i \leq
r$; thus the dominant integral weight $\lambda = \sum_{i=1}^r
m_i \, \omega_i$, where $\omega_i (=\epsilon_1 + \epsilon_2 + \cdots
+ \epsilon_i)$ are the fundamental weights of $\slrpone$. 

Let $\bdf$ be the set of pairs $(\ell,\bds)$ satisfying:
\[ \ell \in \integers_{\geq 0}, \;\; \bds = (\bds(1) \leq \bds(2) \leq \cdots \leq \bds(\ell)) \in \integers_{\geq 0}^\ell. \]
 For $(\ell,\bds) \in \bdf$ and $1 \leq i \leq j \leq r$, define the following element of $U(\slrpone[t])$:
\begin{equation}\label{e:xijcl}
x_{ij}^-(\ell,\bds) =   (x_{ij}^-\tensor t^{\bds(1)}) \, (x_{ij}^-\tensor t^{\bds(2)}) \,\cdots \,  (x_{ij}^-\tensor t^{\bds(\ell)}).
\end{equation}
If $\ell =0$, then the tuple $\bds$ is empty and $x_{ij}^-(\ell,\bds)
=1$. 

Next, we describe the indexing set of the basis of $W(\lambda)$
constructed by Chari-Loktev; this set, denoted $\mathcal{B}^r(\lambda)$, consists of tuples $(\ell_{ij},\bds_{ij})_{1 \leq i \leq j \leq r}$, satisfying the following
properties: (i) $(\ell_{ij},\bds_{ij}) \in \bdf$, and (ii) either $\ell_{ij} =0 $, or $\ell_{ij} >0$ and
\begin{equation} \label{e:clineq}
\bds_{ij}(\ell_{ij}) \leq m_i + \sum_{s=j+1}^r \ell_{i+1,s} - \sum_{s=j}^r \ell_{is}
\end{equation}
for all $1 \leq i \leq j \leq r$ (see equations (2.2) and (2.6) of \cite{cladv2006}). 
Given  $\clpop = (\ell_{ij},\bds_{ij})$ in $\mathcal{B}^r(\lambda)$, define an element $v_\clpop$ of $W(\lambda)$ by:
\begin{equation} \label{e:clbasisvec}
v_\clpop = x_{11}^-(\ell_{11},\bds_{11}) x_{12}^-(\ell_{12},\bds_{12})  x_{22}^-(\ell_{22},\bds_{22})  \cdots x_{rr}^-(\ell_{rr},\bds_{rr}) \; w_{\lambda}
\end{equation}
Note that in the above equation, the terms appear in {\em colexicographic order}, i.e.,  terms involving $x_{ij}^-$ occur to the left of terms involving $x_{i'j'}^-$ if either (i) $j < j'$ or (ii) $j=j'$ and $i < i'$. The vectors $v_\clpop$ form a basis of $W(\lambda)$ as $\clpop$ ranges over $\mathcal{B}^r(\lambda)$ \cite[Theorem 2.1.3]{cladv2006}.

We now establish a bijection between the set $\popset(\lseq)$ of POPs with bounding sequence $\lseq$ and $\mathcal{B}^r(\lambda)$.
This has the further property that if the POP $\pop$ maps to the tuple $\clpop = (\ell_{ij}, \bds_{ij})$, then the vectors $v_\pop$ and $v_\clpop$ are non-zero
scalar multiples of each other. To define the map, let $\pop \in \popset(\lseq)$ be given with underlying pattern $\lseq^1, \lseq^2, \cdots, \lseq^{r+1} = \lseq$ and partition overlay $(\pijiseq)_{1 \leq i \leq j \leq r}$. We then define $\ell_{ij} = \lambda_i^{j+1} - \lambda_i^j$ for all $1 \leq i \leq j \leq r$. Since the partition
$\pijiseq$ fits into the rectangle with sides $(\lambda_i^{j+1} - \lambda_i^j, \lambda_i^{j} - \lambda_{i+1}^{j+1})$, it has at most $\ell_{ij}$ parts. We now let $\bds_{ij}$ be the sequence of parts of $\pijiseq$ arranged in weakly increasing order, padded with an appropriate number of zeros at the beginning so that the length of $\bds_{ij}$ equals $\ell_{ij}$; in other words, if $d_{ij}$ is the number of (non-zero) parts of $\pijiseq$, then $\bds_{ij}(k) = 0$ for $1 \leq k \leq \ell_{ij}-d_{ij}$ and $\bds_{ij}(\ell_{ij}) \geq \cdots \geq \bds_{ij}(\ell_{ij}-d_{ij}+1)$ is the sequence of parts of $\pijiseq$. Finally, define $\clpop(\pop) := (\ell_{ij}, \bds_{ij})_{1 \leq i \leq j \leq r}$. We claim that the map $\pop \mapsto \clpop(\pop)$ is the desired bijection.

We first show that $\clpop(\pop) \in \mathcal{B}^r(\lambda)$. It is clear that $(\ell_{ij},\bds_{ij}) \in \bdf$ for all $i,j$,
so we only need to verify that \eqref{e:clineq} holds.  Now, using the definitions, the right hand side of \eqref{e:clineq} becomes
\begin{align}
m_i + \sum_{s=j+1}^r \ell_{i+1,s} - \sum_{s=j}^r \ell_{is} &=
(\lambda_i^{r+1}-\lambda_{i+1}^{r+1}) + \sum_{s=j+1}^r (\lambda_{i+1}^{s+1} - \lambda_{i+1}^s) -
\sum_{s=j}^r (\lambda_{i}^{s+1} - \lambda_{i}^s)\notag\\  &= \lambda_i^{j} - \lambda_{i+1}^{j+1}. \label{e:clineq_is_gt_se_diff}
\end{align}
But $\bds_{ij}(\ell_{ij})$, being the largest part of $\pijiseq$, is at most $\lambda_i^{j} - \lambda_{i+1}^{j+1}$, thereby establishing \eqref{e:clineq}.
We also observe from definitions that the monomials $x_{ij}^-(\lambda_i^{j+1}-\lambda^j_i,\lambda_i^j-\lambda^{j+1}_{i+1}, \pijiseq)$  of \eqref{e:xijpop} and
$x_{ij}^-(\ell_{ij},\bds_{ij})$ of \eqref{e:xijcl} are scalar multiples of each other; the scaling factor is a product of factorials, since the former involves divided powers while the latter does not.

Now suppose $\clpop  = (\ell_{ij}, \bds_{ij}) \in \mathcal{B}^r(\lambda)$ is given. We construct its inverse $\pop = \pop(\clpop)$ under the above map as follows: (i) The pattern underlying $\pop$ is defined inductively by the relations $\lambda_i^j = \lambda_i^{j+1} - \ell_{ij}$ with $j=r, r-1, \cdots, 1$ and $1\leq i\leq j$. That this defines a pattern follows from equation \eqref{e:clineq_is_gt_se_diff} and the facts that $\ell_{ij} \geq 0$, $\bds_{ij}(\ell_{ij}) \geq 0$.  (ii) The partition $\pijiseq$ in the overlay is simply taken to be sequence $\bds_{ij}$ arranged in weakly decreasing order. That it fits into the appropriate rectangle follows from equation \eqref{e:clineq_is_gt_se_diff} and the fact that $\bds_{ij}$ has length $\ell_{ij}$.

It is clear that the maps $\pop \mapsto \clpop(\pop)$ and $\clpop \mapsto \pop(\clpop)$ are mutual inverses. Further, the remark above concerning the monomials together with equations \eqref{e:clmonom}, \eqref{e:clbasisvec} shows that the vectors $v_\pop$ and $v_{\clpop(\pop)}$ are indeed scalar multiples of each other. This completes the proof of the theorem.
\emyproof

\begin{proposition}\label{p:grade_wt_of_basis}
With notation as above, we have: (i) $v_\pop$ is a homogeneous vector of $W(\lambda)$, with grade equal to the number of boxes in $\pop$, and (ii) the $\lieh$-weight of $v_\pop$ is the weight of (the pattern underlying)~$\pop$.
\end{proposition}
\bmyproof
This is clear from equations \eqref{e:xijpop}, \eqref{e:clmonom}.
\emyproof
Thus we obtain a formula for the graded character of~$W(\lambda)$ (cf. \cite[Prop. 2.1.4]{cladv2006}):
\beq\label{e:grchar}
\boxed{ \textup{char}_tW(\lambda)
=\sum_{\pop}e^{\textup{weight}(\pop)} t^{|\pop|}
}
\eeq
where the sum is over all POPs~$\pop$ with bounding sequence~$\lseq$.

\mysubsection{Representation theoretic proof of
Theorem~\ref{t:main}}\mylabel{ss:rtproof}  We are now ready to give a
representation theoretic proof of the following special case of our main
theorem (Theorem~\ref{t:main}):
\begin{quote}\begin{em}
Let $\lseq: \lambda_1\geq\ldots\geq\lambda_n$ be a non-increasing sequence of integers.    
Let $\mseq\in\mathbb{Z}^n$ such that $\lseq\majgeq\mseq$.   Then there is a unique integral pattern with bounding sequence~$\lseq$,  weight~$\mseq$,  and  area $\frac{1}{2}(||\lseq||^2-||\mseq||^2)$.   Any other pattern (with real entries) with bounding sequence $\lseq$ and weight~$\mseq$ has area strictly less than $\frac{1}{2}(||\lseq||^2-||\mseq||^2)$.
\end{em}\end{quote}

\bmyproof
Subtracting $\lambda_n$ from all entries of a pattern sets up an area and integrality preserving bijection between patterns with bounding sequence $\lseq$ and weight $\mseq$ on the one hand and those with bounding sequence $\lambda_1-\lambda_n\geq\ldots\geq\lambda_{n-1}-\lambda_n\geq0$ and weight $\mseq-(\lambda_n,\ldots,\lambda_n)$ on the other.   Moreover passing from $\lseq$ and $\mseq$ to $\lambda_1-\lambda_n\geq\ldots\lambda_{n-1}-\lambda_n\geq0$
and $\mseq-(\lambda_n,\ldots,\lambda_n)$ does not affect the hypothesis $\lseq\majgeq\mseq$.   We may therefore assume without loss of generality that $\lambda_n=0$.
Let $\lambda, \mu$ denote the weights of $\slrpone$ corresponding to the tuples $\lseq, \mseq$; since $\lseq\majgeq\mseq$, $\mu$ is a weight of the representation $V(\lambda)$.

Now consider the set of integral patterns with bounding sequence $\lseq$ and weight~$\mseq$. If $\pattern$ is one such pattern, then among the POPs which have underlying pattern $\pattern$, there is a unique POP with the largest number of boxes, namely the one in which all partitions in the overlay completely fill up their bounding rectangles. The number of boxes in this POP is clearly the area of $\pattern$. Let $M$ be the maximal area attained among all integral patterns with bounding sequence $\lseq$ and weight~$\mseq$.    Then by Proposition \ref{p:grade_wt_of_basis} and the above discussion,  it follows that~$M$ is maximal such that $\wmlseq_{\mseq}[M]\neq0$ and moreover that the dimension of $\wmlseq_{\mseq}[M]$ equals the number of integral patterns with bounding sequence $\lseq$, weight $\mseq$, and area~$M$.   It now follows from Proposition~\ref{p:key} that the number of such patterns is~$1$ and that $M=\frac{1}{2}(||\lseq||^2-||\mseq||^2)$. This proves the first part of our statement.

Now suppose that we have a pattern with bounding sequence~$\lseq$ and weight~$\mseq$.   Then its trapezoidal area is $\frac{1}{2}(||\lseq||^2-||\mseq||^2)$ (Corollary~\ref{c:traparea}).   So its area is at most this number.  Moreover, if its area equals this number,  then it is integral (Corollary~\ref{c:l:forrtpf}). Thus, the second assertion follows from the first.
\emyproof

\mysection{A bijection between colored partitions and POPs}\mylabel{s:bijection}
\noindent
This section is entirely combinatorial and may be read independently of
the rest of the paper.     Its goal is Theorem~\ref{t:bijection},
which gives a certain bijection between colored partitions of a number
and POPs of a certain kind.      Quite apart from any interest this
bijection may have in its own right,     we use it in the next
section to state the conjectural stability of the Chari-Loktev bases.
The stability property expresses compatibility of the bases with
inclusions of local Weyl modules,   and in order to make sense of this
there must be  in the first place   an identification of the
indexing set of the basis of the included module as a subset of the
indexing set of the basis of the ambient module.      The
combinatorial bijection of this section establishes the desired identification.
%
%
\mysubsection{Breaking up a partition}\mylabel{ss:partbreak}
In this subsection, we describe a procedure to break up a partition
into smaller partitions depending upon some input.   This
can be viewed as a generalization of the construction of Durfee
squares.
We first treat the case when the input is a single integer.   We then treat 
the general case when the input is a non-decreasing sequence of integers.

\mysubsubsection{The case when a single integer is given}\mylabel{sss:pb1}  First suppose that we are given:
\begin{itemize}
\item a partition
$\piseq: \pi_1\geq\pi_2\geq\ldots\ $, and 
\item an integer $c$.
\end{itemize}   It is convenient to put $\pi_0=\infty$.     Consider the function $m\mapsto \pi_m+c-m$ on non-negative integers.   It is decreasing, takes value $\infty$ at $0$,  is non-negative at $c$ if $c$ is non-negative, and is negative for large $m$.   Let $a$ be the largest non-negative integer such that $\pi_a\geq a-c$.  Note that $a\geq c$.

Put $b:=a-c$.    Consider the partitions $\piseqone:\pi_1-b\geq\ldots\geq\pi_a-b$ and $\piseqtwo: \pi_{a+1}\geq \pi_{a+2}\geq\ldots$.   The former has at most $a$ parts;  the latter has largest part at most $b$ (since $\pi_{a+1}<(a+1)-c$ by choice of $a$).   The last assertion can be stated as follows:
\begin{equation}\label{e:pb0}
\textup{$\pi_d\leq b$ for $d>a$}
\end{equation}
It is easily seen that
\begin{equation}\label{e:part1}
|\piseq|=ab+|\piseq^1|+|\piseq^2|
\end{equation}


Consider the association $(c,\piseq)\mapsto(a,b;\piseqone,\piseqtwo)$.   Since $(c,\piseq)$ can be recovered from 
$(a,b;\piseqone,\piseqtwo)$,  the association is one-to-one.   Its
image,  as $c$ varies over all integers and $\piseq$ over all
partitions, consists of all $(a,b;\piseqone,\piseqtwo)$ such that $a$, $b$ are non-negative integers, $\piseqone$~is a partition with at most $a$ parts, and $\piseqtwo$ is a partition with largest part at most $b$.

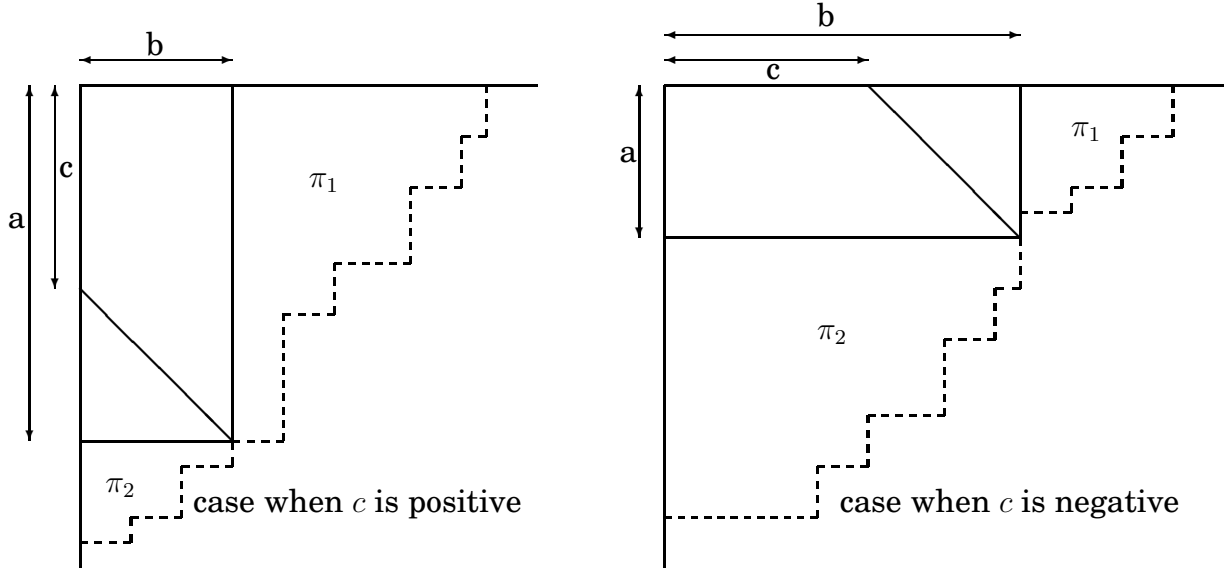
\begin{figure}[h]
\setlength{\unitlength}{1400sp}%
\begingroup\makeatletter\ifx\SetFigFont\undefined%
\gdef\SetFigFont#1#2#3#4#5{%
  \reset@font\fontsize{#1}{#2pt}%
  \fontfamily{#3}\fontseries{#4}\fontshape{#5}%
  \selectfont}%
\fi\endgroup%
\begin{picture}(21378,10083)(706,-9094)
\thicklines
{\color[rgb]{0,0,0}\put(1801,-511){\line( 0,-1){8550}}
}%
{\color[rgb]{0,0,0}\put(1801,-4111){\line( 1,-1){2700}}
}%
{\color[rgb]{0,0,0}\put(4501,-511){\line( 0,-1){6300}}
}%
{\color[rgb]{0,0,0}\put(1801,-6811){\line( 1, 0){2700}}
}%
\thinlines
{\color[rgb]{0,0,0}\put(901,-511){\vector( 0, 1){  0}}
\put(901,-511){\vector( 0,-1){6300}}
}%
{\color[rgb]{0,0,0}\put(1801,-61){\vector(-1, 0){  0}}
\put(1801,-61){\vector( 1, 0){2700}}
}%
{\color[rgb]{0,0,0}\put(1351,-511){\vector( 0, 1){  0}}
\put(1351,-511){\vector( 0,-1){3600}}
}%
\thicklines
{\color[rgb]{0,0,0}\multiput(4501,-6811)(257.14286,0.00000){4}{\line( 1, 0){128.571}}
\multiput(5401,-6811)(0.00000,236.84211){10}{\line( 0, 1){118.421}}
\multiput(5401,-4561)(257.14286,0.00000){4}{\line( 1, 0){128.571}}
\multiput(6301,-4561)(0.00000,257.14286){4}{\line( 0, 1){128.571}}
\multiput(6301,-3661)(245.45455,0.00000){6}{\line( 1, 0){122.727}}
\multiput(7651,-3661)(0.00000,245.45455){6}{\line( 0, 1){122.727}}
\multiput(7651,-2311)(257.14286,0.00000){4}{\line( 1, 0){128.571}}
\multiput(8551,-2311)(0.00000,257.14286){4}{\line( 0, 1){128.571}}
\multiput(8551,-1411)(300.00000,0.00000){2}{\line( 1, 0){150.000}}
\multiput(9001,-1411)(0.00000,257.14286){4}{\line( 0, 1){128.571}}
}%
{\color[rgb]{0,0,0}\put(1801,-511){\line( 1, 0){8100}}
}%
{\color[rgb]{0,0,0}\put(12151,-511){\line( 0,-1){8550}}
}%
{\color[rgb]{0,0,0}\put(18451,-3211){\line(-1, 0){6300}}
}%
{\color[rgb]{0,0,0}\put(18451,-511){\line( 0,-1){2700}}
}%
{\color[rgb]{0,0,0}\put(15751,-511){\line( 1,-1){2700}}
}%
{\color[rgb]{0,0,0}\multiput(18451,-2761)(257.14286,0.00000){4}{\line( 1, 0){128.571}}
\multiput(19351,-2761)(0.00000,300.00000){2}{\line( 0, 1){150.000}}
\multiput(19351,-2311)(257.14286,0.00000){4}{\line( 1, 0){128.571}}
\multiput(20251,-2311)(0.00000,257.14286){4}{\line( 0, 1){128.571}}
\multiput(20251,-1411)(257.14286,0.00000){4}{\line( 1, 0){128.571}}
\multiput(21151,-1411)(0.00000,257.14286){4}{\line( 0, 1){128.571}}
}%
{\color[rgb]{0,0,0}\put(12151,-511){\line( 1, 0){9900}}
}%
\thinlines
{\color[rgb]{0,0,0}\put(11701,-3211){\vector( 0,-1){  0}}
\put(11701,-3211){\vector( 0, 1){2700}}
}%
{\color[rgb]{0,0,0}\put(15751,-61){\vector( 1, 0){  0}}
\put(15751,-61){\vector(-1, 0){3600}}
}%
{\color[rgb]{0,0,0}\put(18451,389){\vector( 1, 0){  0}}
\put(18451,389){\vector(-1, 0){6300}}
}%
\thicklines
{\color[rgb]{0,0,0}\multiput(12151,-8161)(234.78261,0.00000){12}{\line( 1, 0){117.391}}
\multiput(14851,-8161)(0.00000,257.14286){4}{\line( 0, 1){128.571}}
\multiput(14851,-7261)(257.14286,0.00000){4}{\line( 1, 0){128.571}}
\multiput(15751,-7261)(0.00000,257.14286){4}{\line( 0, 1){128.571}}
\multiput(15751,-6361)(245.45455,0.00000){6}{\line( 1, 0){122.727}}
\multiput(17101,-6361)(0.00000,245.45455){6}{\line( 0, 1){122.727}}
\multiput(17101,-5011)(257.14286,0.00000){4}{\line( 1, 0){128.571}}
\multiput(18001,-5011)(0.00000,257.14286){4}{\line( 0, 1){128.571}}
\multiput(18001,-4111)(300.00000,0.00000){2}{\line( 1, 0){150.000}}
\multiput(18451,-4111)(0.00000,257.14286){4}{\line( 0, 1){128.571}}
}%
{\color[rgb]{0,0,0}\multiput(1801,-8611)(257.14286,0.00000){4}{\line( 1, 0){128.571}}
\multiput(2701,-8611)(0.00000,300.00000){2}{\line( 0, 1){150.000}}
\multiput(2701,-8161)(257.14286,0.00000){4}{\line( 1, 0){128.571}}
\multiput(3601,-8161)(0.00000,257.14286){4}{\line( 0, 1){128.571}}
\multiput(3601,-7261)(257.14286,0.00000){4}{\line( 1, 0){128.571}}
\multiput(4501,-7261)(0.00000,300.00000){2}{\line( 0, 1){150.000}}
}%
\put(2971, 29){\makebox(0,0)[lb]{\smash{{\SetFigFont{12}{14.4}{\rmdefault}{\mddefault}{\updefault}{\color[rgb]{0,0,0}b}%
}}}}
\put(1431,-2131){\makebox(0,0)[lb]{\smash{{\SetFigFont{12}{14.4}{\rmdefault}{\mddefault}{\updefault}{\color[rgb]{0,0,0}c}%
}}}}
\put(521,-3031){\makebox(0,0)[lb]{\smash{{\SetFigFont{12}{14.4}{\rmdefault}{\mddefault}{\updefault}{\color[rgb]{0,0,0}a}%
}}}}
\put(5851,-2311){\makebox(0,0)[lb]{\smash{{\SetFigFont{12}{14.4}{\rmdefault}{\mddefault}{\updefault}{\color[rgb]{0,0,0}$\pi_1$}%
}}}}
\put(2251,-7711){\makebox(0,0)[lb]{\smash{{\SetFigFont{12}{14.4}{\rmdefault}{\mddefault}{\updefault}{\color[rgb]{0,0,0}$\pi_2$}%
}}}}
\put(19351,-1411){\makebox(0,0)[lb]{\smash{{\SetFigFont{12}{14.4}{\rmdefault}{\mddefault}{\updefault}{\color[rgb]{0,0,0}$\pi_1$}%
}}}}
\put(14851,539){\makebox(0,0)[lb]{\smash{{\SetFigFont{12}{14.4}{\rmdefault}{\mddefault}{\updefault}{\color[rgb]{0,0,0}b}%
}}}}
\put(13951,-391){\makebox(0,0)[lb]{\smash{{\SetFigFont{12}{14.4}{\rmdefault}{\mddefault}{\updefault}{\color[rgb]{0,0,0}c}%
}}}}
\put(11351,-1861){\makebox(0,0)[lb]{\smash{{\SetFigFont{12}{14.4}{\rmdefault}{\mddefault}{\updefault}{\color[rgb]{0,0,0}a}%
}}}}
\put(14851,-5011){\makebox(0,0)[lb]{\smash{{\SetFigFont{12}{14.4}{\rmdefault}{\mddefault}{\updefault}{\color[rgb]{0,0,0}$\pi_2$}%
}}}}
\put(3801,-8061){\makebox(0,0)[lb]{\smash{{\SetFigFont{12}{14.4}{\rmdefault}{\mddefault}{\updefault}{\color[rgb]{0,0,0}case when $c$ is positive}%
}}}}
\put(15251,-8061){\makebox(0,0)[lb]{\smash{{\SetFigFont{12}{14.4}{\rmdefault}{\mddefault}{\updefault}{\color[rgb]{0,0,0}case when $c$ is negative}%
}}}}
\end{picture}%
\caption{Illustration of the procedure in~\S\ref{sss:pb1}}\label{f:one}
\end{figure}

Figure~\ref{f:one} illustrates the procedure just described.
Note that the $c = 0$ case is the Durfee square construction.
\mysubsubsection{The case when a non-decreasing sequence of integers is given}\mylabel{sss:pbmany}
Now suppose that we are given, for $t\geq2$ an integer, the following:
\begin{itemize}
\item a partition $\piseq: \pi_1\geq\pi_2\geq\ldots\ $, and 
\item a sequence $\cseq: c_1\leq\ldots\leq c_{t-1}$ of integers.
\end{itemize}   
As before,  it is convenient to set $\pi_0=\infty$. 
For $1\leq j\leq t-1$, let $a_j$ be the largest non-negative integer
such that $\pi_{a_j}\geq a_j-c_j$.   Since the $c_j$ are
non-decreasing,  it is clear that the $a_j$ are also non-decreasing:
$a_1\leq \ldots\leq a_{t-1}$.     Set $b_j:=a_j-c_j$.   
\begin{proposition}\mylabel{p:bdecrease}The $b_j$ thus defined are
  non-increasing: $b_1\geq\ldots\geq b_{t-1}$.\end{proposition}    
\bmyproof Fix $j$ such that $1\leq j\leq t-2$ (there is nothing to
prove in case $t=2$).  If $a_{j+1}=a_j$,  then $b_{j+1}=a_{j+1}-c_{j+1}=a_j-c_{j+1}\leq a_j-c_j=b_j$.   
If $a_{j+1}>a_j$,    then,  on the one hand,  $\pi_{a_{j+1}}\leq b_j$ by~(\ref{e:pb0});   and, on the other, $b_{j+1}\leq\pi_{a_{j+1}}$ (by the definitions of $a_{j+1}$ and $b_{j+1}$).\emyproof

\noindent
We define $t$ partitions $\piseqone$, $\piseqtwo$, \ldots, $\piseq^t$ as follows.   Set $a_0=0$, $a_t=\infty$;  $b_0=\infty$, $b_t=0$;    and for $j$, $1\leq j\leq t$:
\begin{equation}\label{e:pij}
\piseqj:\ \pi_{a_{j-1}+1}-b_j\geq \pi_{a_{j-1}+2}-b_j\geq\ldots\geq\pi_{a_j}-b_j
\end{equation}
The above equation can be rewritten as follows:
\begin{equation}\label{e:pij:2}
\textup{$\pi^j_{k-a_{j-1}}:=\pi_k-b_j$ \quad\quad
for $k$ such that $a_{j-1}<k\leq a_j$}
\end{equation}
Note that $\piseqj$ fits into the rectangle
$(a_j-a_{j-1},b_{j-1}-b_j)$ in the sense of~\S\ref{ss:partitions}
(since $\pi_{a_{j-1}+1}\leq b_{j-1}$ by (\ref{e:pb0})).    We have
\begin{equation}\label{e:part:2}
|\piseq|\quad =\quad
|\piseq^1|+\cdots+|\piseq^t|+
\sum_{j=1}^{t-1}(a_j-a_{j-1})b_j\quad =\quad
|\piseq^1|+\cdots+|\piseq^t|+
\sum_{j=1}^{t-1}a_j(b_j-b_{j+1})
\end{equation}

Consider the association $(\cseq,\piseq)\mapsto(\aseq,\bseq;\piseqone,\ldots,\piseq^t)$,  where $\aseq$, $\bseq$ refer respectively to the sequences $a_1\leq\ldots \leq a_{t-1}$ and $b_1\geq\ldots
\geq b_{t-1}$. Since $(\cseq,\piseq)$ can be recovered from
$(\aseq,\bseq;\piseqone,\ldots,\piseq^t)$,  the association is
one-to-one.   Its image, as $\cseq$ varies over all non-decreasing
integer sequences of length $t-1$ and $\piseq$ over all partitions, consists of all $(\aseq,\bseq;\piseqone,\ldots,\piseq^t)$ such that $\aseq$, $\bseq$ are non-negative integer sequences of length $t-1$ with $\aseq$ non-decreasing and $\bseq$ non-increasing,  and for every $j$, $1\leq j\leq t$, $\piseq^j$ is a partition that fits into the rectangle $(a_j-a_{j-1},b_{j-1}-b_j)$.

\begin{figure}[h]
\setlength{\unitlength}{1400sp}%
\begingroup\makeatletter\ifx\SetFigFont\undefined%
\gdef\SetFigFont#1#2#3#4#5{%
  \reset@font\fontsize{#1}{#2pt}%
  \fontfamily{#3}\fontseries{#4}\fontshape{#5}%
  \selectfont}%
\fi\endgroup%
\begin{picture}(20298,13413)(2236,-12244)
\thicklines
{\color[rgb]{0,0,0}\put(3151,389){\line( 1, 0){19350}}
}%
{\color[rgb]{0,0,0}\put(18901,389){\line( 1,-1){1800}}
}%
{\color[rgb]{0,0,0}\put(20701,389){\line( 0,-1){1800}}
}%
{\color[rgb]{0,0,0}\put(15301,389){\line( 0,-1){3150}}
}%
{\color[rgb]{0,0,0}\put(3151,-1411){\line( 1, 0){17550}}
}%
{\color[rgb]{0,0,0}\put(18001,389){\line( 0,-1){3150}}
\put(18001,-2761){\line(-1, 0){14850}}
}%
{\color[rgb]{0,0,0}\put(18001,-2761){\line(-1, 1){3150}}
}%
{\color[rgb]{0,0,0}\put(8551,389){\line( 0,-1){7650}}
\put(8551,-7261){\line(-1, 0){5400}}
}%
{\color[rgb]{0,0,0}\put(8551,-7261){\line( 0,-1){1350}}
\put(8551,-8611){\line(-1, 0){5400}}
}%
{\color[rgb]{0,0,0}\put(8551,-7261){\line(-1, 1){5400}}
}%
{\color[rgb]{0,0,0}\put(8551,-8611){\line(-1, 1){5400}}
}%
{\color[rgb]{0,0,0}\put(5851,389){\line( 0,-1){9900}}
\put(5851,-9511){\line(-1, 0){2700}}
}%
{\color[rgb]{0,0,0}\put(5851,-9511){\line(-1, 1){2700}}
}%
{\color[rgb]{0,0,0}\put(3151,389){\line( 0,-1){12600}}
}%
{\color[rgb]{0,0,0}\multiput(3151,-11761)(300.00000,0.00000){2}{\line( 1, 0){150.000}}
\multiput(3601,-11761)(0.00000,257.14286){4}{\line( 0, 1){128.571}}
\multiput(3601,-10861)(245.45455,0.00000){6}{\line( 1, 0){122.727}}
\multiput(4951,-10861)(0.00000,257.14286){4}{\line( 0, 1){128.571}}
\multiput(4951,-9961)(257.14286,0.00000){4}{\line( 1, 0){128.571}}
\multiput(5851,-9961)(0.00000,300.00000){2}{\line( 0, 1){150.000}}
\multiput(5851,-9511)(257.14286,0.00000){4}{\line( 1, 0){128.571}}
\multiput(6751,-9511)(0.00000,300.00000){2}{\line( 0, 1){150.000}}
\multiput(6751,-9061)(240.00000,0.00000){8}{\line( 1, 0){120.000}}
\multiput(8551,-9061)(0.00000,300.00000){2}{\line( 0, 1){150.000}}
\multiput(8551,-8611)(0.00000,245.45455){6}{\line( 0, 1){122.727}}
\multiput(8551,-7261)(257.14286,0.00000){4}{\line( 1, 0){128.571}}
\multiput(9451,-7261)(0.00000,257.14286){4}{\line( 0, 1){128.571}}
\multiput(9451,-6361)(240.00000,0.00000){8}{\line( 1, 0){120.000}}
\multiput(11251,-6361)(0.00000,300.00000){2}{\line( 0, 1){150.000}}
\multiput(11251,-5911)(257.14286,0.00000){4}{\line( 1, 0){128.571}}
\multiput(12151,-5911)(0.00000,300.00000){2}{\line( 0, 1){150.000}}
}%
{\color[rgb]{0,0,0}\put(15301,-2761){\line(-1, 1){3150}}
}%
{\color[rgb]{0,0,0}\multiput(13501,-4111)(257.14286,0.00000){4}{\line( 1, 0){128.571}}
\multiput(14401,-4111)(0.00000,257.14286){4}{\line( 0, 1){128.571}}
\multiput(14401,-3211)(257.14286,0.00000){4}{\line( 1, 0){128.571}}
\multiput(15301,-3211)(0.00000,300.00000){2}{\line( 0, 1){150.000}}
\multiput(15301,-2761)(232.25806,0.00000){16}{\line( 1, 0){116.129}}
\multiput(18901,-2761)(0.00000,300.00000){2}{\line( 0, 1){150.000}}
\multiput(18901,-2311)(300.00000,0.00000){2}{\line( 1, 0){150.000}}
\multiput(19351,-2311)(0.00000,300.00000){2}{\line( 0, 1){150.000}}
\multiput(19351,-1861)(245.45455,0.00000){6}{\line( 1, 0){122.727}}
\multiput(20701,-1861)(0.00000,300.00000){2}{\line( 0, 1){150.000}}
\multiput(20701,-1411)(257.14286,0.00000){4}{\line( 1, 0){128.571}}
\multiput(21601,-1411)(0.00000,257.14286){4}{\line( 0, 1){128.571}}
\multiput(21601,-511)(300.00000,0.00000){2}{\line( 1, 0){150.000}}
\multiput(22051,-511)(0.00000,257.14286){4}{\line( 0, 1){128.571}}
}%
{\color[rgb]{0,0,0}\multiput(12331,-5371)(296.53868,296.53868){4}{\line( 1, 1){145.384}}
}%
\thinlines
{\color[rgb]{0,0,0}\put(5851,839){\vector( 1, 0){0}}
\put(3151,839){\vector(-1, 0){0}}
\multiput(3160,839)(9.00000,0.00000){299}{\makebox(1.5875,11.1125){\tiny.}}
}%
{\color[rgb]{0,0,0}\put(8551,839){\vector( 1, 0){0}}
\put(5851,839){\vector(-1, 0){0}}
\multiput(5860,839)(9.00000,0.00000){299}{\makebox(1.5875,11.1125){\tiny.}}
}%
{\color[rgb]{0,0,0}\put(18001,839){\vector( 1, 0){0}}
\put(15301,839){\vector(-1, 0){0}}
\multiput(15310,839)(9.00000,0.00000){299}{\makebox(1.5875,11.1125){\tiny.}}
}%
{\color[rgb]{0,0,0}\put(20701,839){\vector( 1, 0){0}}
\put(18001,839){\vector(-1, 0){0}}
\multiput(18010,839)(9.00000,0.00000){299}{\makebox(1.5875,11.1125){\tiny.}}
}%
{\color[rgb]{0,0,0}\put(2701,-1411){\vector( 0,-1){0}}
\put(2701,389){\vector( 0, 1){0}}
\multiput(2701,380)(0.00000,-9.00000){199}{\makebox(1.5875,11.1125){\tiny.}}
}%
{\color[rgb]{0,0,0}\put(2701,-2761){\vector( 0,-1){0}}
\put(2701,-1411){\vector( 0, 1){0}}
\multiput(2701,-1420)(0.00000,-9.00000){149}{\makebox(1.5875,11.1125){\tiny.}}
}%
{\color[rgb]{0,0,0}\put(2701,-8611){\vector( 0,-1){0}}
\put(2701,-7261){\vector( 0, 1){0}}
\multiput(2701,-7270)(0.00000,-9.00000){149}{\makebox(1.5875,11.1125){\tiny.}}
}%
{\color[rgb]{0,0,0}\put(2701,-9511){\vector( 0,-1){0}}
\put(2701,-8611){\vector( 0, 1){0}}
\multiput(2701,-8620)(0.00000,-9.00000){99}{\makebox(1.5875,11.1125){\tiny.}}
}%
\put(651,-3031){\makebox(0,0)[lb]{\smash{{\SetFigFont{12}{14.4}{\rmdefault}{\mddefault}{\updefault}{\color[rgb]{0,0,0}$\boxed{a'_3=0}$}%
}}}}
\put(2131,-511){\makebox(0,0)[lb]{\smash{{\SetFigFont{12}{14.4}{\rmdefault}{\mddefault}{\updefault}{\color[rgb]{0,0,0}$a'_1$}%
}}}}
\put(2131,-2041){\makebox(0,0)[lb]{\smash{{\SetFigFont{12}{14.4}{\rmdefault}{\mddefault}{\updefault}{\color[rgb]{0,0,0}$a'_2$}%
}}}}
\put(4411,1019){\makebox(0,0)[lb]{\smash{{\SetFigFont{12}{14.4}{\rmdefault}{\mddefault}{\updefault}{\color[rgb]{0,0,0}$b'_t$}%
}}}}
\put(7111,1019){\makebox(0,0)[lb]{\smash{{\SetFigFont{12}{14.4}{\rmdefault}{\mddefault}{\updefault}{\color[rgb]{0,0,0}$b'_{t-1}$}%
}}}}
\put(8111,1319){\makebox(0,0)[lb]{\smash{{\SetFigFont{12}{14.4}{\rmdefault}{\mddefault}{\updefault}{\color[rgb]{0,0,0}$\boxed{b'_{t-2}=0}$}%
}}}}
\put(16561,1019){\makebox(0,0)[lb]{\smash{{\SetFigFont{12}{14.4}{\rmdefault}{\mddefault}{\updefault}{\color[rgb]{0,0,0}$b'_3$}%
}}}}
\put(19261,1019){\makebox(0,0)[lb]{\smash{{\SetFigFont{12}{14.4}{\rmdefault}{\mddefault}{\updefault}{\color[rgb]{0,0,0}$b'_2$}%
}}}}
\put(1541,-7891){\makebox(0,0)[lb]{\smash{{\SetFigFont{12}{14.4}{\rmdefault}{\mddefault}{\updefault}{\color[rgb]{0,0,0}$a'_{t-2}$}%
}}}}
\put(1541,-9151){\makebox(0,0)[lb]{\smash{{\SetFigFont{12}{14.4}{\rmdefault}{\mddefault}{\updefault}{\color[rgb]{0,0,0}$a'_{t-1}$}%
}}}}
\put(3781,-10141){\makebox(0,0)[lb]{\smash{{\SetFigFont{12}{14.4}{\rmdefault}{\mddefault}{\updefault}{\color[rgb]{0,0,0}$\piseq^t$}%
}}}}
\put(6101,-9200){\makebox(0,0)[lb]{\smash{{\SetFigFont{12}{14.4}{\rmdefault}{\mddefault}{\updefault}{\color[rgb]{0,0,0}$\piseq^{t-1}$}%
}}}}
\put(18361,-1951){\makebox(0,0)[lb]{\smash{{\SetFigFont{12}{14.4}{\rmdefault}{\mddefault}{\updefault}{\color[rgb]{0,0,0}$\piseq^2$}%
}}}}
\put(21151,-331){\makebox(0,0)[lb]{\smash{{\SetFigFont{12}{14.4}{\rmdefault}{\mddefault}{\updefault}{\color[rgb]{0,0,0}$\piseq^1$}%
}}}}
\put(18611, 19){\makebox(0,0)[lb]{\smash{{\SetFigFont{12}{14.4}{\rmdefault}{\mddefault}{\updefault}{\color[rgb]{0,0,0}$c_1$}%
}}}}
\put(14471,19){\makebox(0,0)[lb]{\smash{{\SetFigFont{12}{14.4}{\rmdefault}{\mddefault}{\updefault}{\color[rgb]{0,0,0}$c_2$}%
}}}}
\put(11771,19){\makebox(0,0)[lb]{\smash{{\SetFigFont{12}{14.4}{\rmdefault}{\mddefault}{\updefault}{\color[rgb]{0,0,0}$c_3$}%
}}}}
\put(3331,-1861){\makebox(0,0)[lb]{\smash{{\SetFigFont{12}{14.4}{\rmdefault}{\mddefault}{\updefault}{\color[rgb]{0,0,0}$c_{t-3}$}%
}}}}
\put(3331,-3211){\makebox(0,0)[lb]{\smash{{\SetFigFont{12}{14.4}{\rmdefault}{\mddefault}{\updefault}{\color[rgb]{0,0,0}$c_{t-2}$}%
}}}}
\put(3331,-6811){\makebox(0,0)[lb]{\smash{{\SetFigFont{12}{14.4}{\rmdefault}{\mddefault}{\updefault}{\color[rgb]{0,0,0}$c_{t-1}$}%
}}}}
\put(11541,-6711){\makebox(0,0)[lb]{\smash{{\SetFigFont{12}{14.4}{\rmdefault}{\mddefault}{\updefault}{\color[rgb]{0,0,0}In
this instance, $c_1$,
$c_2$, $c_3$ are negative,}%
}}}}
\put(11541,-7411){\makebox(0,0)[lb]{\smash{{\SetFigFont{12}{14.4}{\rmdefault}{\mddefault}{\updefault}{\color[rgb]{0,0,0}and
$c_{t-3}$, $c_{t-2}$, $c_{t-1}$ are positive}%
}}}}
\put(11541,-9451){\makebox(0,0)[lb]{\smash{{\SetFigFont{12}{14.4}{\rmdefault}{\mddefault}{\updefault}{\color[rgb]{0,0,0}
$a'_j:=a_j-a_{j-1}$ and $b'_j:=b_{j-1}-b_j$}%
}}}}
\put(11541,-10151){\makebox(0,0)[lb]{\smash{{\SetFigFont{12}{14.4}{\rmdefault}{\mddefault}{\updefault}{\color[rgb]{0,0,0}
so that $\piseq^j$ fits into $(a'_j,b'_j)$}%
}}}}
\put(12541,-5051){\makebox(0,0)[lb]{\smash{{\SetFigFont{12}{14.4}{\rmdefault}{\mddefault}{\updefault}{\color[rgb]{0,0,0}$\iddots$}%
}}}}
\put(8741,-8451){\makebox(0,0)[lb]{\smash{{\SetFigFont{12}{14.4}{\rmdefault}{\mddefault}{\updefault}{\color[rgb]{0,0,0}$\boxed{\textup{$\piseq^{t-2}$
is empty}}$}%
}}}}
\put(15741,-3451){\makebox(0,0)[lb]{\smash{{\SetFigFont{12}{14.4}{\rmdefault}{\mddefault}{\updefault}{\color[rgb]{0,0,0}$\boxed{\textup{$\piseq^{3}$
is empty}}$}%
}}}}
\end{picture}%
\caption{Illustration of the procedure in~\S\ref{sss:pbmany}}\label{f:two}
\end{figure}
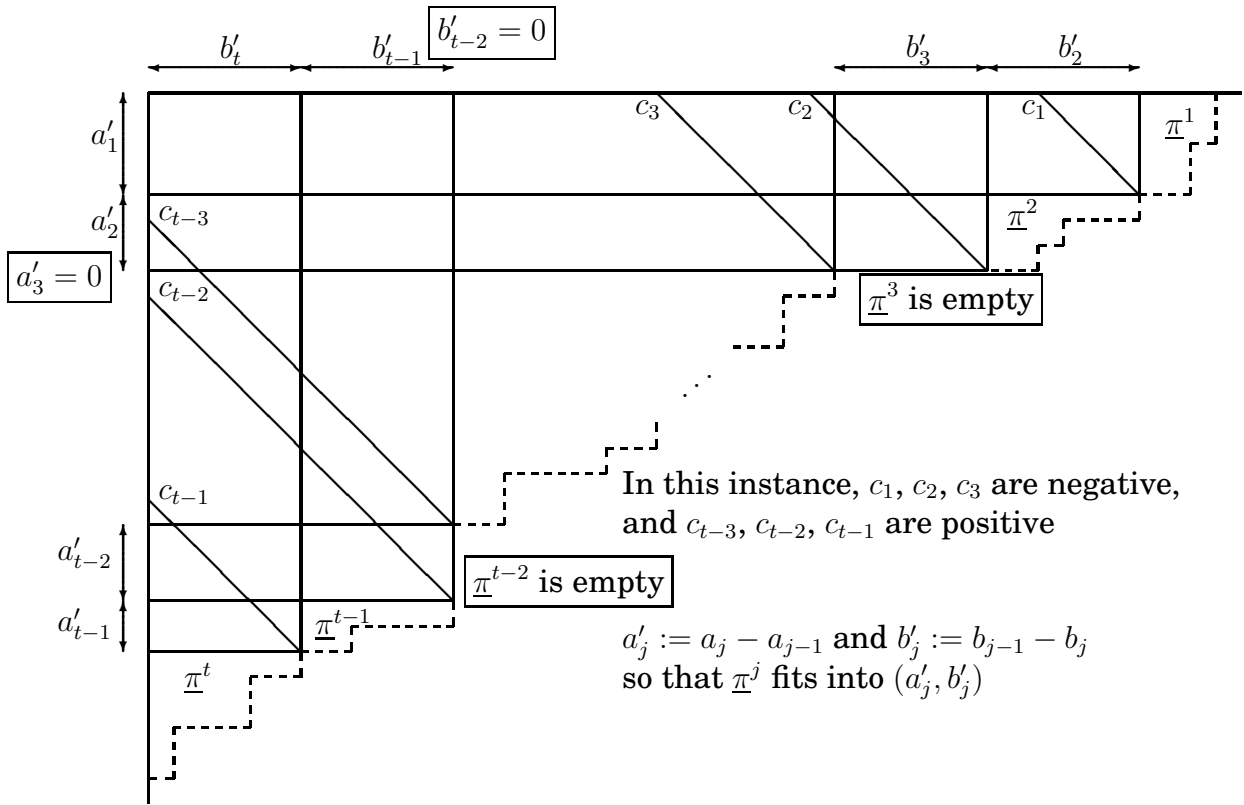

The picture in Figure~\ref{f:two} describes the procedure just described.
\mysubsection{Nearly interlacing
  sequences with approximate partition overlays}\mylabel{ss:ilseq} 
Let $s\geq1$ be an integer and $\etaseq$: $\eta_1$,
\ldots, $\eta_{s+1}$ with $\eta_2\geq\ldots\geq\eta_s$ an integer sequence.
An integer sequence $\etatwoseq$: $\etatwo_1$, \ldots, $\etatwo_s$ with
$\etatwo_2\geq\ldots\geq\etatwo_s$ is said to {\em nearly interlace}
$\etaseq$ if either $s=1$ (in which case no further condition is
imposed) or $s\geq2$ and
\begin{equation}\label{e:d:nearil}\textup{$\eta_1\geq\etatwo_1$,\quad\quad
$\etatwo_s\geq\eta_{s+1}$,\quad\quad  and 
$\etatwo_2\geq\ldots\geq\etatwo_{s-1}$
  interlaces $\eta_2\geq\cdots\geq\eta_s$.}\end{equation}
The 
following is a pictorial depiction of this definition (where $x\longrightarrow y$ means 
$x\geq y$,   and $\times$ indicates the absence of any relation):
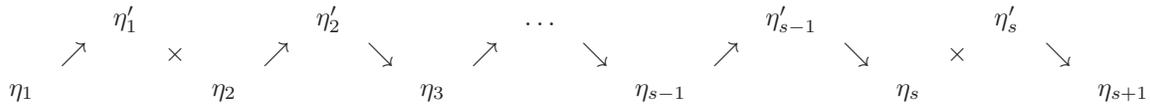
\begin{equation*}
\begin{array}{ccccccccccccccccccccc}
&&\etatwo_1&&&&\etatwo_2&&&&\ldots&&&&\etatwo_{s-1}&&&&\etatwo_s\\
&\nearrow&&\times&&\nearrow&&\searrow&&\nearrow&&\searrow&&\nearrow&&\searrow&&\times&&\searrow\\
\eta_1&&&&\eta_2&&&&\eta_3&&&&\eta_{s-1}&&&&\eta_s&&&&\eta_{s+1}\\
\end{array}
\end{equation*}
In case $s=1$,  the definition imposes no further constraint on
$\etatwoseq$ and so the pictorial depiction is:
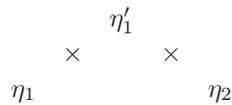
\begin{equation*}
\begin{array}{ccccc}
&&\etatwo_1\\
&\times&&\times\\
\eta_1&&&&\eta_2\\
\end{array}
\end{equation*}
We define the {\em proper trapezoidal area\/} of the nearly 
interlacing sequences $\etaseq$, $\etatwoseq$ as above by:
\begin{equation}
\label{e:d:proptrap}
\proptrap{\etaseq,\etatwoseq}:=\sum_{1\leq i<j\leq s}(\eta_i-\etatwo_i)(\etatwo_j-\eta_{j+1}) 
\end{equation}
Observe that this is non-negative in general and further that it is
zero if $s=1$.

Given sequences $\etaseq$, $\etatwoseq$ as above that nearly
interlace,  a 
sequence $\piseq^1$, \ldots, $\piseq^s$ of partitions is said to {\em
approximately  overlay\/} $\etaseq$, $\etatwoseq$,   if
either $s=1$ (in which case no further condition is imposed) or
$s\geq2$ and
\begin{equation}\label{e:d:overlay}
\left\{\begin{array}{l}
\textup{$\piseq^1$ has at most $\eta_1-\etatwo_1$ parts,  
$\piseq^s$ has largest part at most $\etatwo_s-\eta_{s+1}$, and}\\
\textup{for $j=2,\ldots, s-1$,   the partition $\piseq^j$ fits into the 
  rectangle $(\eta_j-\etatwo_j,\etatwo_j-\eta_{j+1})$.}\\
\end{array}\right.\end{equation}

\mysubsubsection{Producing nearly interlacing
  sequences with approximate partition overlays}\mylabel{sss:ilseq} 
Fix an integer $s\geq1$ and an integer sequence $\etaseq$: $\eta_1$,
\ldots, $\eta_{s+1}$ with $\eta_2\geq\ldots\geq\eta_s$.     
Suppose that we are given:   
\begin{itemize}
\item a partition $\piseq:\pi_1\geq \pi_2\geq\ldots\ $,  and
 \item an integer $\mu$.
\end{itemize}  
Our goal first of all in this subsection is to associate to the 
data $(\mu, \piseq)$ an integer  sequence $\etaseq'$ that nearly
interlaces $\etaseq$ together
with an approximate partition overlay on $\etaseq$, $\etaseq'$.
The map is denoted 
by~$\bijxieta$.   We then investigate
the nature of $\bijxieta$ (Lemma~\ref{l:bijection}).

Put $c_1:=\mu-\eta_2$,  
\ldots, $c_{s-1}:=\mu-\eta_s$.
Then $c_1\leq\ldots\leq c_{s-1}$.    We apply the procedure
of~\S\ref{sss:pbmany} with $t=s$ and $\cseq$ as above to obtain
$(\aseq, \bseq,\piseq^1,\ldots,\piseq^s)$.     In case~$s=1$,   we
take $\aseq$ and $\bseq$  to be empty and set $\piseq^1:=\piseq$.
In case $s=2$,  the procedure of~\S\ref{sss:pbmany} reduces to that
of~\S\ref{sss:pb1}.


We now define the sequence $\etatwoseq$. 
In case $s=1$,  set
\begin{equation}\label{e:etatwosone}
\etatwo_1:=\eta_1+\eta_2-\mu\quad\quad\quad\quad\quad\quad\quad\textup{(case  $s=1$)}
\end{equation}     Now suppose $s\geq2$.
As before,  it is convenient to set $a_0=b_s=0$ and $a_s=b_0=\infty$.    Define
\begin{equation}\label{e:etatwo}
\etatwo_j:=\eta_j-(a_j-a_{j-1}) \quad\textup{for $j=1,\ldots,s-1$} \quad\textup{and}\quad\quad
\etatwo_s:=\eta_{s+1}+b_{s-1}\quad\quad\textup{(case $s\geq2$)}
\end{equation}
\begin{proposition}\mylabel{p:etatwo} Suppose that $s\geq 2$.   Then:
\begin{enumerate}
\item For $j=2,\ldots,s$,   we have $\etatwo_j=\eta_{j+1}+(b_{j-1}-b_j)$.
\item $\etaseq$ and $\etatwoseq$ are nearly interlaced.
\item The sequence $\piseq^1$, \ldots, $\piseq^s$ of partitions
  approximately overlays $\etaseq$,
  $\etatwoseq$.
\item
$(\eta_1+\cdots+\eta_{s+1})-(\etatwo_1+\cdots+\etatwo_s)=\mu$
\item
$\proptrap{\etaseq,\etatwoseq}|+(|\piseq^1|+\cdots+|\piseq^s|) \ =\
|\piseq| $
\item $\sum_{i=1}^j(\eta_i-\etatwo_i) =a_j$ for $0\leq j<s$ and
  $\sum_{i=j+1}^s(\etatwo_i-\eta_{i+1})=b_j$ for $1\leq j\leq s$.
\end{enumerate}
\end{proposition}
\bmyproof  For $j=s$,~(1) is just the definition of~$\etatwo_s$.
 Fix $j$ in the range $2$, \ldots, $s-1$.     We have:
\begin{align*}
\etatwo_j & = \eta_j-(a_j-a_{j-1})  & \textup{definition of $\etatwo_j$}\\
& = \eta_j-((c_j+b_j)-(c_{j-1}+b_{j-1})) & \textup{definition of $\bseq$}\\
& = \eta_j-((\mu-\eta_{j+1}+b_j)-(\mu-\eta_j+b_{j-1})) & \textup{definition of $\cseq$}\\
& = \eta_{j+1}+(b_{j-1}-b_j)
\end{align*}
This proves (1).    For (2),  we observe:
\begin{itemize}
\item For $1\leq j\leq s-1$, we have $\eta_j\geq\eta_j-(a_j-a_{j-1})=\etatwo_j$  since $\underline{a}$ is a non-decreasing sequence.
\item For $2\leq j\leq s$, we have
  $\etatwo_j=\eta_{j+1}+(b_{j-1}-b_j)\geq\eta_{j+1}$  by (1) and the
  fact that $\underline{b}$ is a non-increasing sequence (Proposition~\ref{p:bdecrease}).
\end{itemize}
Assertion (3) follows since $\piseq^j$ fits into $(a_j-a_{j-1},b_{j-1}-b_j)$ by construction.

Assertion (4) is just the definition
(\ref{e:etatwosone}) in case $s=1$.     Now suppose $s\geq 2$. 
By the definition (\ref{e:etatwo}) of $\etatwo_j$, we have 
\begin{align*}
\etatwo_1+\cdots+\etatwo_s & = (\eta_1-(a_1-a_0))\ +\ \cdots\ +\ (\eta_{s-1}-(a_{s-1}-a_{s-2}))\quad+\quad (\eta_{s+1}+b_{s-1})\\
& = \eta_1+\cdots+\eta_{s-1}+\eta_{s+1}-a_{s-1}+b_{s-1}   \\
& = \eta_1+\cdots+\eta_{s-1}+\eta_{s+1}-c_{s-1} \quad\textup{(since $b_{s-1}=a_{s-1}-c_{s-1}$ by definition)}\\
& = \eta_1+\cdots+\eta_{s-1}+\eta_s+\eta_{s+1}-\mu \quad\textup{(since
  $c_{s-1}=\mu-\eta_s$ by definition)
}
\end{align*}
For (5),  first rewrite the definition (\ref{e:d:proptrap}) to get:
\begin{equation*}
\proptrap{\etaseq,\etatwoseq}\quad=\quad\sum_{1\leq i\leq s-1}
(\eta_i-\etatwo_i)\left(\sum_{i+1\leq j\leq 
    s}(\etatwo_j-\eta_{j+1})\right) 
\end{equation*}
Substituting from (\ref{e:etatwo}) and item (1) 
into the right hand side above,  we get 
\begin{equation*}
\proptrap{\etaseq,\etatwoseq}\quad=\quad 
\sum_{1\leq i\leq s-1}
(a_i-a_{i-1})\left(\sum_{i+1\leq j\leq 
    s}(b_{j-1}-b_{j})\right) 
 \quad =\quad 
\sum_{1\leq i\leq s-1}
(a_i-a_{i-1})b_i 
\end{equation*}
Assertion (5) now follows from (\ref{e:part:2}).

For (6),   using (\ref{e:etatwo}) we obtain
$\sum_{i=1}^j(\eta_i-\etatwo_i)=\sum_{i=1}^j(a_i-a_{i-1})=a_j$.
Similarly,  using (1) we obtain
$\sum_{i=j+1}^s(\etatwo_i-\eta_{i+1})=\sum_{i=j+1}^s(b_{i-1}-b_i)=b_j$.
\emyproof



\blemma\mylabel{l:bijection}
Fix an integer $s\geq1$ and an integer sequence $\etaseq$: $\eta_1$,
\ldots, $\eta_{s+1}$  with $\eta_2\geq\ldots\geq\eta_{s}$.
Let $\bijxieta$ denote the association described above:
\begin{itemize}
\item
from the set of all pairs $(\mu,
\piseq)$,  where $\mu$ is an integer and $\piseq$  a partition
\item to the set of all tuples $(\etatwoseq,
  \piseq^1,\ldots,\piseq^s)$,  where $\etatwo$: $\etatwo_1$, \ldots,
  $\etatwo_s$ is an integer sequence nearly interlacing $\etaseq$,
  and $\piseq^1$, \ldots, $\piseq^s$  a sequence of partitions
  approximately overlaying $\etaseq$, $\etatwoseq$
\end{itemize}
The association $\bijxieta$ is a bijection.     More precisely,  the
association $\bijxipeta$ in the other direction to be defined below
(in the course of the proof of this lemma) is the two-sided inverse of $\bijxieta$.
\elemma 
\bmyproof
We define $\bijxipeta$.    Let $\etatwoseq$ and $\piseq^1$, \ldots,
$\piseq^s$ with the specified properties be given.       The image of 
$(\etatwoseq, \piseq^1, \ldots, \piseq^s)$ under $\bijxipeta$ is
defined to be $(\mustar, \piseqstar)$ where $\mustar$ and $\piseqstar$
are as defined below.   Set 
\begin{equation}\label{e:d:mustar}\mustar:=(\eta_1+\cdots+\eta_{s+1})-(\etatwo_1+\cdots+\etatwo_s)
\end{equation}
To define $\piseqstar$,   we take $k$ to be a positive integer and define
$\pistar_k$.     For $j$, $0\leq j\leq s-1$,   set
$\astar_j:=\sum_{i=1}^j(\eta_i-\etatwo_i)$.  Put $\astar_s:=\infty$.    We have
$0=\astar_0\leq \astar_1\leq \ldots\leq \astar_{s-1}<\astar_s=\infty$.   Thus there exists
unique $j$, $1\leq j\leq s$, such that $\astar_{j-1}<k\leq \astar_j$.      Set
\begin{equation}\label{e:pistar}
\textup{$\pistar_k:=\pi^j_{k-\astar_{j-1}}+\bstar_j$ \quad\quad 
where $\bstar_j:=\sum_{i=j+1}^s(\etatwo_i-\eta_{i+1})$}
\end{equation}
Since $\etatwo_i\geq\eta_{i+1}$ for all $i\geq2$,  it is clear that
$\bstar_j$ and hence also $\pistar_k$ is non-negative.

Let us verify that $\pistar_k\geq\pistar_{k+1}$ for all $k$.   If
$\astar_{j-1}<k<\astar_j$,   then $\astar_{j-1}<k+1\leq\astar_j$,  so
that, from (\ref{e:pistar}),
$\pistar_k-\pistar_{k+1}=\pi^j_{k-\astar_{j-1}}-\pi^j_{k+1-\astar_{j-1}}\geq
0$ (since $\piseq^j$ is a partition).    Now suppose that
$k=\astar_j$.   Then, from (\ref{e:pistar}),    
\begin{equation*}
\pistar_k-\pistar_{k+1}=(\pi^j_{k-\astar_{j-1}}+\bstar_j)-(\pi^{j+1}_1+\bstar_{j+1})=\pi^j_{\astar_j-\astar_{j-1}}+((\etatwo_{j+1}-\eta_{j+2})-\pi^{j+1}_1)\geq0
\end{equation*}
where the last inequality holds since $\piseq^{j+1}$ has largest part
at most $\etatwo_{j+1}-\eta_{j+2}$ by hypothesis.
This proves that $\pistarseq$ is a partition and finishes the definition of $\bijxipeta$.

We now verify that $\bijxipeta\circ\bijxieta$ is the identity.
Suppose we first apply $\bijxieta$ to $(\mu,\piseq)$ to get
$(\etatwoseq,\piseq^1,\ldots,\piseq^s)$ to which in turn we
apply~$\bijxipeta$ to get $(\mustar,\piseqstar)$.     From (4) of
Proposition~\ref{p:etatwo} and (\ref{e:d:mustar}) it follows that
$\mustar=\mu$.     It follows from the definitions of $\astarseq$ and
$\bstarseq$ above and (6) of Proposition~\ref{p:etatwo} that
$\astarseq=\aseq$ and $\bstarseq=\bseq$.    It now follows from the
definitions (\ref{e:pij:2}) and (\ref{e:pistar}) respectively of
$\piseq^j$ and $\pistarseq$ that $\pistarseq=\piseq$.

Finally we verify that $\bijxieta\circ\bijxipeta$ is the identity.
Let $(\mustar,\piseqstar)$ be the result of application of
$\bijxipeta$ to $(\etatwoseq, \piseq^1,\ldots,\piseq^s)$.    To
calculate the action of $\bijxieta$ on $(\mustar,\piseqstar)$,  we
must compute $\cseq$, $\aseq$, and $\bseq$ as in~\S\ref{sss:pbmany}.
From the definition of~$\cseq$ at the beginning of this subsection and
those of $\mustar$, $\astarseq$, and $\bstarseq$ above,  we have:
\begin{equation}\label{e:cj}
c_j=\mustar-\eta_{j+1}=(\sum_{i=1}^j(\eta_i-\etatwo_{i}))-(\sum_{i=j+1}^{s}(\etatwo_i-\eta_{i+1}))=\astar_j-\bstar_j
\end{equation}
We claim that $\astarseq=\aseq$.    Assuming this claim,  it follows
from (\ref{e:cj}) and the definition of $\bseq$ in~\S\ref{sss:pbmany}
that $\bseq=\bstarseq$.     From (\ref{e:etatwo}) and the definition
of $\astarseq$ above,  it follows that the nearly interlacing sequence
part of the image of $\bijxieta$ is $\etatwoseq$.
From (\ref{e:pij:2}), (\ref{e:pistar}), and the equalities
$\aseq=\astarseq$, $\bseq=\bstarseq$,  it now follows that the image
under $\bijxieta$ of $(\mustar,\piseqstar)$ equals
$(\etatwoseq,\piseq^1,\ldots,\piseq^s)$.

It remains only to prove the claim above that $\aseq=\astarseq$, or in
other words that for $j$, $1\leq j\leq s-1$,   $\astar_j$ is the
largest integer such that $\pistar_{\astar_j}\geq \astar_j-c_j$.
From (\ref{e:pistar}) and (\ref{e:cj}),    we have 
\[\pistar_{\astar_j}\ =\ \pi^j_{\astar_j-\astar_{j-1}}+\bstar_j\ \geq
\ \bstar_j\ =\ \astar_j-c_j\]
We now show that $\pistar_{\astar_j+1}<\astar_j+1-c_j$.   Fix $\ell$,
$j+1\leq\ell\leq s$, such that
$\astar_j=\astar_{l-1}<\astar_j+1\leq\astar_\ell$.  From
(\ref{e:pistar}) and (\ref{e:cj}), we get
\begin{align*}
1+\astar_j-c_j-\pistar_{\astar_j+1}& = 1
+(\astar_j-c_j)-\pi^\ell_1-\bstar_\ell\\  &= 1 +(\bstar_j-\bstar_\ell)
-\pi^\ell_1\\
&= 1+ (\etatwo_{j+1}-\eta_{j+2})+\cdots+(\etatwo_{\ell-1}-\eta_\ell)+((\etatwo_{\ell}-\eta_{\ell+1})-\pi^\ell_1)
\end{align*}
Since $(\etatwo_i-\eta_{i+1})\geq0$ for $2\leq i\leq s$ and the
largest part $\pi^\ell_1$ of $\piseq^\ell$ is at most
$(\etatwo_\ell-\eta_{\ell+1})$,   the right hand side in the last
line of the above
display is positive.
\emyproof
\mysubsection{Near patterns with approximate partition overlays}\mylabel{ss:npattern}
Fix an integer $r\geq1$.    Suppose that, for every $j$, $1\leq j\leq r+1$, we
have an integer sequence $\lseqj$: $\lj_1$, \ldots, $\lj_{j}$ of
length $j$ with $\lj_2\geq\ldots\geq\lj_{j-1}$.    We say that this
collection of sequences forms a {\em near pattern\/}
if $\lseq^{j}$ nearly interlaces $\lseq^{j+1}$ for every $j$, $1\leq
j\leq r$.     The last sequence $\lseqrpone$ is called the {\em bounding
  sequence\/} of the near pattern.    
The 
following is a pictorial depiction of this definition for $r=4$ (where $x\longrightarrow y$ means 
$x\geq y$,   and $\times$ indicates the absence of any relation):
\begin{equation*}
\begin{array}{ccccccccccccccccccccc}
&&&&&&&&&&\lambda_1^1\\
&&&&&&&&&\times&&\times\\
&&&&&&&&\lambda_1^2&&&&\lambda_2^2\\
&&&&&&&\nearrow&&\times&&\times&&\searrow\\
&&&&&&\lambda^3_1&&&&\lambda^3_2&&&&\lambda^3_3\\
&&&&&\nearrow&&\times&&\nearrow&&\searrow&&\times&&\searrow\\
&&&&\lambda^4_1&&&&\lambda^4_2&&&&\lambda^4_3&&&&\lambda^4_4\\
&&&\nearrow&&\times&&\nearrow&&\searrow&&\nearrow&&\searrow&&\times&&\searrow\\
&&\lambda_1^5&&&&\lambda_2^5&&&&\lambda_3^5&&&&\lambda_4^5&&&&\lambda^5_5\\
\end{array}
\end{equation*}
The {\em proper trapezoidal area\/} of a near pattern $\pattern=\{\ljseq\st
1\leq j\leq r+1\}$ is defined by:
\begin{equation}\label{e:d:proptrap:np}
\proptrap{\pattern}\ :=\ \sum_{j=2}^r \proptrap{\ljponeseq,\ljseq}
\ =\ \sum_{j=2}^r \sum_{1\leq i<h\leq j} (\lambda^{j+1}_i-\lambda^j_i)(\lambda^j_h-\lambda^{j+1}_{h+1})
\end{equation}
The {\em weight\/} of a near pattern~$\pattern$ as above is the tuple
$(\mu_1,\ldots,\mu_{r+1})$,    where 
\[
\textup{$\mu_{j+1}:=\sum_{i=1}^{j+1}\lambda^{j+1}_i -
  \sum_{i=1}^j\lambda^j_i$ \quad for $1\leq j\leq r$
  \quad\quad\quad\quad and \quad $\mu_1:=\lambda^1_1$}
\]

Let $\pattern=\{\lseqj\st 1\leq j\leq r+1\}$ be a near pattern.    Suppose that
we are given partitions $\pijiseq$,   for $1\leq j\leq r$ and
$1\leq i\leq j$.    We say that this collection of partitions {\em approximately
  overlays\/} the near pattern~$\pattern$ if:
\begin{itemize}
\item for $2\leq j\leq r$,   the partition $\pijoneseq$ has 
  at most $\ljpone_1-\lj_1$ parts
\item for $2\leq j\leq r$,   the partition $\pijseq^j$ has
  largest part at most $\lj_j-\ljpone_{j+1}$ 
\item for $3\leq j\leq r$ and $2\leq i\leq j-1$,  the partition
  $\pijiseq$ fits into the rectangle $(\ljpone_{i}-\lj_i,\lj_i-\ljpone_{i+1})$
\end{itemize}
The above conditions can also be expressed by saying that for
every $j$, $1\leq j\leq r$,  the sequence $\pijiseq$, $1\leq i\leq j$,
of partitions approximately
overlays the nearly interlacing sequences $\lseqjpone$, $\lseqj$ in
the sense of (\ref{e:d:overlay}).

The {\em number of boxes\/} in an approximate partition overlay as
above of a near pattern  is defined to be $\sum_{1\leq i\leq j\leq
  r}|\pijiseq|$.   The terminology is justified by thinking of the
partitions in terms of their shapes.

We sometimes use the the term {\em approximately overlaid near
  pattern\/}, {\em AONP\/} for short,  for
a near pattern with an approximate partition overlay.

\mysubsubsection{A bijection on AONPs}\mylabel{sss:bijnp}
Fix an integer $s\geq1$ and an integer sequence $\lseqrpone$:
$\lrpone_1$, \ldots, $\lrpone_{r+1}$ with
$\lrpone_2\geq\ldots\geq\lrpone_r$.   
Let  $\nplrpone$ denote the set of all AONPs with bounding sequence $\lrponeseq$.    
Let $\intpart$ denote the set of
tuples $(\mseq; \pioneseq, \ldots, \pirseq)$, where
$\mseq$: $\mutwo$, \ldots, $\murpone$ is a sequence of $r$ integers and
$\pioneseq$, \ldots, $\pirseq$ is a sequence of $r$ partitions (the
reason for the indexing of $\mu_j$ starting with $2$ will become
clear presently).

Given
an element $(\mseq; \pioneseq, \ldots, \pirseq)$ in $\intpart$,
set  $(\lrseq,\pirseq^1,\ldots,\pirseq^r):=
\bijxilrpone(\murpone,\pirseq)$,  where $\bijxilrpone$ is as defined
in~\S\ref{sss:ilseq}.
By (2) of~Proposition~\ref{p:etatwo},  it follows that $\lrseq$ is
such that $\lr_2\geq\ldots\geq\lr_{r-1}$.    We may thus inductively
define:
\begin{equation}\label{e:ljpijseq}
(\ljseq,\pijseq^1,\ldots,\pijseq^j):=\bijxiljpone(\mu_{j+1},\pijseq) 
\end{equation}
From (2) and (3) of Proposition~\ref{p:etatwo},  it follows that the
sequences $\ljseq$ ($1\leq j\leq r+1$) and partitions $\pijseq^i$
($1\leq j\leq r$, $1\leq i\leq j$) form an AONP with bounding sequence $\lrponeseq$.    Thus we have defined a
map from $\intpart$ to $\nplrpone$,  which too we denote by
$\bijxilrpone$ by abuse of notation.
\blemma\mylabel{l:bijnp}
The map $\bijxilrpone$ from $\intpart$ to $\nplrpone$ just defined is
a bijection.   
\elemma
\bmyproof
We construct a two
sided inverse. 
For an element
$\{\ljseq; \pijseq^i\st 1\leq j\leq r, 1\leq i\leq j\}$ of
$\nplrpone$,    set 
\[
(\mustar_{j+1},\pijseq^\star) :=
\bijxipljpone(\ljseq,\pijseq^1,\ldots,\pijseq^j)
\]
where $\bijxipljpone$ is as defined in the proof of
Lemma~\ref{l:bijection}.    Since $\bijxipljpone$ is the two sided inverse
of $\bijxiljpone$ (by the assertion of that lemma),   it follows that the
map \[\{\ljseq; \pijseq^j\}\mapsto
(\mustar_{2},\ldots,\mustar_{r+1};\pioneseq^\star, \ldots,
\pirseq^\star)\] is the required two sided inverse.    By abusing
notation again,  we denote it by $\bijxiplrpone$.
\emyproof
\bprop\mylabel{p:wtbox}  
The underlying near pattern~$\pattern$ in the image  of 
$(\mu_2,\ldots,\mu_{r+1};\pioneseq,\ldots,\pirseq)$ under
$\bijxilrpone$  has bounding sequence  $\lrponeseq$   and weight
$(\mu_1,\mu_2,\ldots,\mu_{r+1})$ where
$\mu_1:=(\sum_{j=1}^{r+1}\lrpone_j)-(\sum_{j=2}^{r+1}\mu_j)$.  The
number~$n$ of boxes in the approximate partition overlay satisfies:
\begin{equation}\label{e:box}
\proptrap{\pattern}+n\ =\ |\pioneseq|+\cdots+|\pirseq|
\end{equation}
\eprop
\bmyproof   The first assertion is immediate from the definition.  From the definition of~$\bijxilrpone$ and
Proposition~\ref{p:etatwo}~(4) and~(5),  we obtain, for $j$ such that $1\leq j\leq r$,
$(\sum_{i=1}^{j+1}\lambda^{j+1}_i)
-(\sum_{i=1}^j\lambda_i^j)=\mu_{j+1}$ and
$\proptrap{\lseq^{j+1},\lseq^j}+\sum_{i=1}^j|\pijseq^i|=|\pijseq|$.
Adding the first set of equations, we get
$\lambda_1^1=(\sum_{j=1}^{r+1}\lrpone_j)-(\sum_{j=2}^{r+1}\mu_j)=\mu_1$.
Adding the second, we get (\ref{e:box}).
\emyproof

\mysubsection{Near patterns to patterns: the shift map}\mylabel{ss:shift}
\noindent
Let an integer $k$ be fixed.     Given a sequence $\etaseq$: $\eta_1$,
\ldots, $\eta_{s+1}$,   where ($s\geq0$ is an integer)  we denote by
$\letaseq$ the sequence $\leta_1$, \ldots, $\leta_{s+1}$,  where
\begin{equation}\label{e:d:leta}
\left\{
\begin{array}{llr}
\leta_1:=\eta_1+2k \quad\quad\quad &\leta_j:=\eta_j+k\quad  \textup{for $2\leq j\leq
  s$}\quad\quad \textup{and}\ \leta_{s+1}:=\eta_{s+1}  \quad\quad&
\textup{(if $s\geq1$)}\\
&\leta_1:=\eta_1+k  & \textup{(if $s=0$)}
\end{array}
\right.
\end{equation}
We refer to $\letaseq$ as the {\em shift by $k$\/} (or just {\em
  shift\/} if $k$ is clear from the context) of $\etaseq$.
The suppression of the dependence on~$k$ in the notation $\letaseq$
should cause no confusion.

The shift $\letaseq$, $\letatwoseq$ of a pair 
$\letaseq$, $\letatwoseq$ of nearly interlacing sequences is also nearly
interlacing.    A sequence of partitions approximately overlays $\etaseq$,
$\etatwoseq$ if and only if it approximately overlays $\letaseq$, $\letatwoseq$.

The {\em shift of a near pattern\/} consists of the shifts of
the constituent sequences of the pattern.     It is also a near
pattern.    A collection of partitions approximately overlays a near pattern if and
only if it approximately overlays the shifted pattern.

To shift a pair of nearly interlacing sequences or a near
pattern with an approximate partition overlay, we just shift the
constituent sequences.   The partitions in the
overlay stay as they are.

Shift preserves proper trapezoidal area (of a pair of nearly interlacing
sequences or of a near pattern).      If a near pattern has weight $(\mu_1,\ldots,\mu_{s+1})$,
its shift by $k$ has weight $(\mu_1+k,\ldots,
\mu_{s+1}+k)$.  Positive shifts of interlacing sequences (respectively
patterns) continue to be interlacing (respectively patterns).
Given a pair of nearly interlacing sequences (respectively a near pattern),  its
shift by $k$ for $k\gg0$ is a pair of  interlacing sequences
(respectively a pattern).
\bprop\mylabel{p:shift}
Let an integer $s\geq1$ and an integer sequence $\etaseq$: $\eta_1$,
\ldots, $\eta_{s+1}$ with $\eta_2\geq\ldots\geq\eta_s$ be fixed.
For $\mu$ an integer and $\piseq$ a partition,  if
$\bijxi_{\etaseq}(\mu,\piseq) =(\etatwoseq, \piseq^1, \ldots, \piseq^s)$,
then $\bijxi_{\letaseq}(\mu+k,\piseq)=(\letatwoseq,\piseq^1,\ldots,\piseq^s)$.
\eprop
\bmyproof
Indeed, the sequences $\cseq$, $\aseq$, and $\bseq$
involved in the calculation of
$\bijxi_{\etaseq}(\mu,\piseq)$ (see~\S\ref{sss:pbmany}) are the same as the corresponding ones
involved in the calculation of $\bijxi_{\letaseq}(\mu+k,\piseq)$.
The result now follows from the definition (\ref{e:pij}) of
$\piseq^j$ and (\ref{e:etatwosone}), (\ref{e:etatwo}) of $\etatwoseq$.
\emyproof
\blemma\mylabel{l:niltoil}
Let $\etaseq$:  $\eta_1$, \ldots, $\eta_{s+1}$ and $\etatwoseq$:
$\etatwo_1$, \ldots, $\etatwo_s$ be
integer sequences (for some $s\geq 1$),  let $\mseq=(\mu_1,\ldots,
\mu_{s+1})$ be a tuple of integers, and $\piseq^1$, \ldots, $\piseq^s$
a sequence of partitions.   Assume that:
\begin{enumerate}
\item $\etaseq$ is non-increasing:  $\eta_1\geq\ldots\geq\eta_{s+1}$
\item $\etaseq$, $\etatwoseq$ nearly interlace and $\piseq^1$, \ldots,
  $\piseq^s$ is an approximate partition overlay on $\etaseq$, $\etatwoseq$
\item $\etaseq\majgeq\mseq$
\item $(\eta_1+\cdots+\eta_{s+1})-(\etatwo_1+\cdots+\etatwo_s)=\mu_{s+1}$
\end{enumerate}
Let $k$ be an integer and $\letaseq$, $\letatwoseq$ be the shifts by $k$
of $\etaseq$, $\etatwoseq$.   If $k\geq\proptrap{\etaseq,
  \etatwoseq}+|\piseq^1|+|\piseq^s|$,   then 
\begin{enumerate}
\item[(a)] $\letaseq$, $\letatwoseq$ interlace (which in particular 
  implies that $\letatwoseq$ is non-increasing:
  $\letatwo_1\geq\ldots\geq\letatwo_s$) 
\item[(b)] $\piseq^1$, \ldots, $\piseq^s$ overlays $\letaseq$,
  $\letatwoseq$ 
(i.e., $\piseq^j$ fits  into the rectangle $(\leta_j-\letatwo_j,
\letatwo_j-\leta_{j+1})$ for $1\leq j\leq s$) 
\item[(c)] $\letatwoseq\majgeq(\mu_1+k, \ldots,\mu_s+k)$
\end{enumerate}
\bmyproof  
For (b),  since $\piseq^1$, \ldots, $\piseq^s$ approximately overlays
$\letaseq$, $\letatwoseq$,  it is enough to show:
\begin{equation}\label{e:l:niltoil:0}
\textup{$\piseq^1$ has largest part at most $\letatwo_1-\leta_2$\quad\quad and \quad\quad
$\piseq^s$ has at most $\leta_s-\letatwo_s$ parts}\end{equation}
These assertions imply in particular that $\letatwo_1\geq\leta_2$ and
$\leta_s\geq\letatwo_s$,  from which (a) follows (since $\letaseq$,
$\letatwoseq$ nearly interlace).

To prove (\ref{e:l:niltoil:0}), we first show that $\letatwo_1-\leta_2\geq\pi^1_1$.   Putting 
$\letatwo_1=\etatwo_1+2k$ and $\leta_2=\eta_2+k$,    we see that the 
desired inequality is equivalent to 
\begin{equation}\label{e:l:niltoil:1}(\eta_1-\eta_2)+k\geq (\eta_1-\etatwo_1)+\pi^1_1 
\end{equation}
We consider two cases.   First suppose that $\etatwo_j-\eta_{j+1}\geq 1$ for 
some $2\leq j\leq s$.   Then, since 
\[k\ \geq\ \proptrap{\etaseq,\etatwoseq}+|\piseq^1|\ \geq\
(\eta_1-\etatwo_1)(\etatwo_j-\eta_{j+1})+\pi^1_1\ \geq\
(\eta_1-\etatwo_1)+\pi^1_1\]
and $(\eta_1-\eta_2)\geq0$ by hypothesis~(1),   we obtain 
(\ref{e:l:niltoil:1}).    In the second case,  we have 
$\etatwo_j=\eta_{j+1}$ for all $j\geq2$.    Then, by hypothesis~(4),
we get $\eta_1-\etatwo_1=\mu_{s+1}-\eta_2$.  Substituting this into 
(\ref{e:l:niltoil:1}),  we get 
\begin{equation}\label{e:l:niltoil:2}
(\eta_1-\mu_{s+1}) + k \geq \pi^1_1 
\end{equation}
But since $k\geq|\piseq^1|\geq \pi^1_1$ and $\eta_1\geq\mu_{s+1}$ by 
hypothesis (3),  we obtain (\ref{e:l:niltoil:2}).

The proof of the latter half of (\ref{e:l:niltoil:0}),  
namely that $\letatwo_s-\leta_{s+1}\geq p$  (where $p$ denotes
the number of parts in $\piseq^s$), is analogous to that of the former
half above. Putting $\leta_2=\eta_s+k$
$\letatwo_s=\etatwo_s$,    we see that the 
desired inequality is equivalent to 
\begin{equation}\label{e:l:niltoil:3}(\eta_s-\eta_{s+1})+k\geq (\etatwo_s-\eta_{s+1})+p
\end{equation}
We consider two cases.   First suppose that $\eta_j-\etatwo_{j}\geq 1$ for 
some $1\leq j<s$.   Then, since 
\[k\ \geq\ \proptrap{\etaseq,\etatwoseq}+|\piseq^s|\ \geq\
(\eta_j-\etatwo_j)(\etatwo_s-\eta_{s+1})+p\ \geq\
(\etatwo_s-\eta_{s+1})+p\]
and $(\eta_s-\eta_{s+1})\geq0$ by hypothesis~(1),   we obtain 
(\ref{e:l:niltoil:3}).    In the second case,  we have 
$\eta_j=\etatwo_{j}$ for all $j<s$.    Then, by hypothesis~(4),
we get $\etatwo_s-\eta_{s+1}=\eta_s-\mu_{s+1}$.  Substituting this into 
(\ref{e:l:niltoil:3}),  we get 
\begin{equation}\label{e:l:niltoil:4:old}
(\mu_{s+1}-\eta_{s+1}) + k \geq p 
\end{equation}
But since $k\geq|\piseq^s|\geq p$ and $\mu_{s+1}\geq\eta_{s+1}$ by 
hypothesis (3),  we obtain (\ref{e:l:niltoil:4:old}).

We now turn to (c).   That
$\letatwo_1+\cdots+\letatwo_s=(\mu_1+k)+\cdots+(\mu_s+k)$ follows from
the definition of $\letatwoseq$, hypothesis (4), and the fact implied
by hypothesis~(3) that
$\eta_1+\cdots+\eta_{s+1}=\mu_1+\cdots+\mu_{s+1}$.   It follows from
(a) that $\letatwo_1\geq\ldots\geq\letatwo_s$.    It only remains to
show that for any $j<s$ and $1\leq i_1<\ldots<i_j\leq s$ we have
\begin{equation}\label{e:l:niltoil:5}
\letatwo_1+\cdots+\letatwo_j\geq (\mu_{i_1}+k)+\cdots+(\mu_{i_j}+k)
\end{equation}
Substituting the definitions $\letatwo_1=\etatwo_1+2k$ and $\letatwo_i=\etatwo_i+k$
for $1<i\leq j$,    we may rewrite (\ref{e:l:niltoil:5}) equivalently
as
\begin{equation}\label{e:l:niltoil:6}
(\eta_1+\cdots+\eta_j)-(\mu_{i_1}+\cdots+\mu_{i_j}) + k\geq (\eta_1-\etatwo_1)+\cdots+(\eta_j-\etatwo_j)
\end{equation}
We consider two cases.   First suppose that $\etatwo_i-\eta_{i+1}\geq 1$ for 
some $j<i\leq s$.   Then, since 
\[k\ \geq\ \proptrap{\etaseq,\etatwoseq}\ \geq\
\left((\eta_1-\etatwo_1)+\cdots+(\eta_j-\etatwo_j)\right)(\etatwo_i-\eta_{i+1})\ \geq\
(\eta_1-\etatwo_1)+\cdots+(\eta_j-\etatwo_j)
\]
and $(\eta_1+\cdots+\eta_j)-(\mu_{i_1}+\cdots+\mu_{i_j})\geq0$ by hypothesis~(3),   we obtain 
(\ref{e:l:niltoil:6}).    In the second case,  we have 
$\etatwo_i=\eta_{i+1}$ for all $i>j$.    Then, by hypothesis~(4),
we get $(\eta_1-\etatwo_1)+\cdots+(\eta_j-\etatwo_j)=\mu_{s+1}-\eta_{j+1}$.  Substituting this into 
(\ref{e:l:niltoil:6}),  we get 
\begin{equation}\label{e:l:niltoil:4}
(\eta_1+\cdots+\eta_{j+1}) - (\mu_{i_1}+\cdots+\mu_{i_j}+\mu_{s+1})+ k \geq 0
\end{equation}
But  since $k\geq0$ and $(\eta_1+\cdots+\eta_{j+1})\geq
(\mu_{i_1}+\cdots+\mu_{i_j}+\mu_{s+1})$ by
hypothesis (3),  we are done.
\emyproof
\elemma
\begin{corollary}\mylabel{c:aonptopop}
Let $\lseq:  \lambda_1\geq\ldots\geq\lambda_{r+1}$ be a non-increasing
integer sequence (where $r\geq 1$ is an integer).      Let
$\mseq\in\mathbb{Z}^{r+1}$ be such that $\lseq\majgeq\mseq$.    Let
$\pattern$ be a near pattern with bounding sequence $\lseq$ and weight
$\mseq$.    Let $\{\pijiseq\st 1\leq i\leq j\leq r\}$ be an
approximate partition overlay on~$\pattern$.      Then the shift by
$\proptrap{\pattern}+\sum_{1\leq j\leq r}(|\pijseq^1|+|\pijseq^j|)$ of the AONP
$(\pattern, \{\pijiseq\})$ is a POP.
\end{corollary}
\bmyproof   We proceed by induction on~$r$.     Let $\ljseq$, $1\leq
j\leq r+1$,  be the constituent sequences of $\pattern$.   We have
$\lseq=\lrponeseq$.   Consider the shift $\lpattern$ of $\pattern$ by
$k=\proptrap{\lrponeseq, \lrseq}+|\pirseq^1|+|\pirseq^r|$.    By the lemma,   we obtain:
\begin{enumerate}
\item[(a)] $\llrponeseq$, $\llrseq$ interlace (which in particular 
  implies that $\llrseq$ is non-increasing:
  $\llr_1\geq\ldots\geq\llr_s$) 
\item[(b)] $\pirseq^1$, \ldots, $\pirseq^r$ overlays $\llrponeseq$,
  $\llrseq$ 
\item[(c)] $\llrseq\majgeq(\mu_1+k, \ldots,\mu_r+k)$
\end{enumerate}
In particular,  this gives a proof in the base case $r=1$ of the
induction.       

Now suppose $r\geq2$.     Let $\patone$ denote the
pattern obtained from $\pattern$ by omitting its last row.      We may
apply the induction hypothesis to $\llrseq$,
$(\mu_1+k,\ldots,\mu_r+k)$, and the AONP $(\lpattern_1, \{\pijiseq\st 1\leq i\leq
j\leq r-1\})$,  to conclude that the shift by
$\proptrap{\lpattern_1}+\sum_{1\leq
  j<r}(|\pijseq^1|+|\pijseq^j|)$    
of this AONP is a POP.

Note the following:
\begin{itemize}
\item $\square^{\textup{prop}}$ is preserved under shifts
\item $\proptrap{\pattern}=\proptrap{\lrponeseq,\lrseq}+\proptrap{\pattern_1}$
\item positive shifts of $\llrponeseq$, $\llrseq$ do
not affect (a) and (b)
\item shift by $\proptrap{\pattern_1}+\sum_{1\leq
    j<r}(|\pijseq^1|+|\pijseq^j|)$ of $\lpattern_1$ (respectively $\llrponeseq$, $\llrseq$) equals shift by
  $\proptrap{\pattern}+\sum_{1\leq j\leq r}(|\pijseq^1|+|\pijseq^j|)$
  of~$\pattern_1$ (respectively $\lrponeseq$, $\lrseq$)
\end{itemize}
The result follows.
\emyproof
\mysubsection{Bijection between $r$-colored partitions and
  POPs}\mylabel{ss:bijection}
We are at last ready to state and prove the desired theorem.    
Fix integers $r\geq1$ and $d\geq0$.  Let $\rpartd$ denote the set of all 
$r$-colored partitions of~$d$. 

Fix $\lseqrpone=\lseq: \lambda_1\geq\ldots\geq\lambda_{r+1}$
a non-increasing sequence of integers,   and  
$\mseq=(\mu_1,\ldots,\mu_{r+1})\in \mathbb{Z}^{r+1}$ such that
$\lseq\majgeq\mseq$.      Let  us denote a general AONP with bounding
sequence of length $r+1$ by $(\pattern, \{\pijiseq\st1\leq i\leq j\leq
r\})$,  where $\pattern$ denotes the underlying near pattern and $\pijiseq$
the partitions in the approximate overlay.
Let $\nlmd$ denote the set of all AONPs where 
$\pattern$ has bounding sequence $\lrponeseq$,
  weight $(\mu_1,\ldots,\mu_{r+1})$, and satisfies the following condition:
\begin{equation}\label{e:boxcon}\proptrap{\pattern}+\sum_{1\leq i\leq
    j\leq r}|\pijiseq|=d
\end{equation}

Let $\philm$ be the map from $\rpartd$ to the set of AONPs given by
\[ \textup{$\pioneseq$, \ldots, $\pirseq$ \quad $\mapsto$ \quad
  $\bijxilrpone(\mu_2,\ldots,\mu_{r+1}; \pioneseq,\ldots,\pirseq)$}\]
where $\bijxilrpone$ is as defined in~\S\ref{sss:bijnp}.  
\bprop\mylabel{p:philm}
The map $\philm$ is a bijection from $\rpartd$ to $\nlmd$.
\eprop
\bmyproof
 It follows
from Proposition~\ref{p:wtbox} that the image of $\philm$ lies
in~$\nlmd$.   Since $\bijxilrpone$ is a bijection
(Lemma~\ref{l:bijnp}),  it follows that $\philm$ is an injection.
We now show that is also onto~$\nlmd$.    Given an element of $\nlmd$,
its image under $\bijxiplrpone$ maps to that element under
$\bijxilrpone$ (see the proof of Lemma~\ref{l:bijnp}),   so it is of
the form $(\mu_2,\ldots,\mu_{r+1};\pioneseq,\ldots,\pirseq)$,  where 
$\pioneseq$, \ldots, $\pirseq$ is an $r$-colored partition of~$d$
(Proposition~\ref{p:wtbox}).
\emyproof

For an integer $k$,  let $\shiftk$ denote the ``shift by $k$''
operator (\S\ref{ss:shift}).    Let $\nlmkd$ (respectively $\plmkd$) denote the set of
all AONPs (respectively POPs) where 
$\pattern$ has bounding sequence $\llseq$  (the shift by $k$ of 
$\lseq$),
  weight $(\mu_1+k,\ldots,\mu_{r+1}+k)$, and (\ref{e:boxcon}) is
  satisfied. 
\begin{theorem}\mylabel{t:bijection}   For $k\geq d$, the composition
  $\shiftk\circ\philm$ defines a bijection from $\rpartd$ to $\plmkd$.
\end{theorem}
\bmyproof   The operator $\shiftk$ is evidently a bijection from
$\nlmd$ to $\nlmkd$  (see~\S\ref{ss:shift}).   By
Corollary~\ref{c:aonptopop},  its image lies in $\plmkd$,    so
$\nlmkd=\plmkd$.   The result follows from Proposition~\ref{p:philm}.
\emyproof

\mysubsection{The complementation involution $\comp$}\mylabel{ss:comp}
For $\pop=(\pattern, \pijiseq\st 1\leq i\leq j\leq r)$ a POP with
$\etaseq^j$ being the $j^\textup{th}$ row of the pattern~$\pattern$,
let $\barpop$ denote the POP $(\pattern, \pijciseq\st 1\leq i\leq j\leq 
r)$, where $\pijciseq$ denotes the complement of $\pijiseq$ in the 
rectangle $(\eta^{j+1}_i-\eta^j_i,\eta^j_i-\eta^{j+1}_{i+1})$
(see~\S\ref{sss:comppart}).   For $\pop$ in~$\plmkd$, we have
\begin{equation*}
\sum_{1\leq i\leq j\leq r}|\pijciseq|
= \area{\pattern}-\sum_{1\leq i\leq j\leq r}|\pijiseq|
=\area{\pattern}+\proptrap{\pattern}-d 
=\traparea{\pattern}-d 
\end{equation*} 
The association $\comp:\pop\mapsto\barpop$ is evidently reversible.   The
above calculation shows that~$\comp$ defines a bijection from $\plmkd$ onto $\plmksd$,  where
$\plmksd$ denotes the set of POPs with bounding sequence~$\llseq$,
weight $(\mu_1+k,\ldots,\mu_{r+1}+k)$,  and depth $d$ (see~\S\ref{ss:pop}
for the definition of depth).

Precomposing the bijection of Theorem~\ref{t:bijection} with $\comp$
continues to be a bijection:
\bcor\mylabel{c:t:bijection}
For $k\geq d$, the composite map $\comp\circ\shiftk\circ\philm$ is a bijection from
$\rpartd$ to $\plmksd$.
\ecor
\bprop\mylabel{p:compat}
The bijections of Theorem~\ref{t:bijection} and
Corollary~\ref{c:t:bijection} are compatible.   More precisely,    for
$j\geq0$,   we have:
\begin{gather} 
\label{e:compat:1}\shiftjpk\circ\philm=\shiftj\circ\phishift\quad\quad\quad \textup{in the theorem}
\\
\label{e:compat:2}\comp\circ\shiftjpk\circ\philm=\comp\circ\shiftj\circ\phishift \quad\quad\quad\textup{in the corollary}
\end{gather}
\eprop
\bmyproof   The left hand side of (\ref{e:compat:2})  may be written
as $\comp\circ\shiftj\circ(\shiftk\circ\philm)$.   Now, by Proposition~\ref{p:shift},
$\shiftk\circ\philm=\phishift$.   The proof of~(\ref{e:compat:1}) is similar.
\emyproof
\bcor\mylabel{c:c:bijection}
For $k\geq d$ and $j\geq0$,  the map $\comp\circ\shiftj\circ\comp:
\plmksd\to\plmkpjsd$ is a bijection.
\ecor
\bmyproof
$\comp\circ\shiftj\circ\comp$ is equal to the composition of two
bijections:   the inverse of $\comp\circ\shiftk\circ\philm$ followed
by $\comp\circ\shiftjpk\circ\philm$.
\emyproof

\mysection{The stability conjecture}\mylabel{s:stabnew}
\noindent
In this section we state the conjecture about stability of
Chari-Loktev bases under a chain of inclusions of local Weyl modules (for $\lieg=\slrpone$).    The
conjecture has been proved in~\cite{rrv:stab} in the case $r=1$.%
\footnote{The conjecture has since been proved by one of us~\cite{br:stab}.}
We begin by recalling details about the chain of inclusions.     The
theorem of the previous section (\S\ref{s:bijection}) gives us an
identification of the indexing set of the basis for an
included module as a subset of the indexing set of the basis of the
larger module.     It then makes sense to ask whether the bases are
well behaved with respect to inclusions.   The stability conjecture
says that this is so in the stable range.

%
\mysubsection{The set up}\mylabel{ss:s1:setup}
Let $\lambda$ be a dominant integral weight and $\theta$ the highest root of~$\lieg=\slrpone$.
We identity $\lambda$ and $\theta$ with $(r+1)$-tuples of integers $\lseq$ and $\tseq$ respectively as in~\S\ref{ss:wttup}:  $\lseq$: $\lambda_1\geq\cdots\lambda_{r}\geq\lambda_{r+1}=0$ and
$\tseq=(2,1,\ldots,1,0)$.  
%
%
%
\mysubsubsection{Chains of inclusions}\mylabel{sss:chn:s1}
As recalled in~\S\ref{ss:lwmdem},  local Weyl modules are Demazure modules.  It follows from this---see~\cite[Lemma~8]{fladv2007} or, for more detail, \cite[\S5.1.2]{br:thesis}---that there is a chain of $\curralg$-module inclusions as follows, each of which is defined uniquely up to scaling:
\[ W(\lambda)\hookrightarrow W(\lambda+\theta)\hookrightarrow 
W(\lambda+2\theta)\hookrightarrow\ldots\]
%
%
The Chari-Loktev basis for
$W(\lambda+k\theta)$ is indexed by the set $\pkl$ of POPs with
bounding sequence $\lseq+k\tseq$  (\S\ref{ss:clbase}).    Mirroring
the above chain of inclusions of local Weyl modules,  we have a chain
of  inclusions of these indexing sets:
\[\pnotl\hookrightarrow\ponel\hookrightarrow\ptwol\hookrightarrow\ldots\]
Indeed, the map $\comp\shiftj\comp$,  where $\comp$ is the
complementation (\S\ref{ss:comp}) and $\shiftj$ is the shift by $j$
(\S\ref{ss:shift}),  defines an injection $\pil\hookrightarrow\pipjl$.
These injections are compatible as $i$ and $j$ vary over the non-negative integers (since
$\comp$ is an involution and $\shiftj\shift^{j'}=\shift^{j+j'}$).
For $\pop \in \pil$, we let $\pop^j$ denote its image under the inclusion
$\pil\hookrightarrow\pipjl$ ($i,j \geq 0$).

\mysubsubsection{Scaling of generators}\mylabel{sss:scale:s1}
For a POP~$\pop$,   let $\vcl{\pop}$ and $\rhocl{\pop}$ denote
respectively the Chari-Loktev basis element and the Chari-Loktev
monomial corresponding to~$\pop$.     Note that the Chari-Loktev basis
depends upon the choice of the generator (which we may scale).     Fix
a generator $w_\lambda$ of $\wml$.  For every $k\geq1$,  let the generator $w_{\lambda+k\theta}$ of $W(\lambda+k\theta)$ be scaled so that:
\beq\label{e:scale}
w_{\lambda}\mapsto\vcl{\popnotk}=\rhocl{\popnotk}w_{\lambda+k\theta}
\quad\quad\quad\textup{under the inclusion $W(\lambda)\hookrightarrow W(\lambda+k\theta)$}
\eeq
where $\popnot$ denotes  the unique element of $\pnotl$ of weight
$\lseq$ (corresponding to the generator of $W(\lambda)$).

\mysubsection{Stability conjecture}\mylabel{ss:stab:first}
It is now natural to ask whether,  for all $k\geq0$, 
and for all $\pop\in\pkl$, we have
\beq\label{e:stab:0}
\vcl{\pop}\mapsto \pm \vcl{\pop^j} \quad \text{ for all } j\geq0,\quad
\quad\textup{under the inclusions $W(\lambda+k\theta)\hookrightarrow W(\lambda+(j+k)\theta)$} 
\eeq
Simple instances (see~\cite[Example~1, \S3.3]{rrv:stab}) show
that~(\ref{e:stab:0}) is too much to expect in general.   We do
however conjecture that it holds when $k$ is in the ``stable range'':

\begin{conjecture}[{\scshape Stability} of Chari-Loktev bases]  
With notation as above,  let $\pop$ be a POP in~$\pkl$.    Let $\mseq$
be the weight of~$\pop$ and $d$ its depth.    Note that
$\mseq\majleq\lseq+k\tseq$ and in particular
$\sum_{i=1}^{r+1}\mu_i=(\sum_{i=1}^{r+1}\lambda_i)+k(r+1)$.     The
assertion (\ref{e:stab:0}) holds if $k\geq \ell+d$, where $\ell$ be the least non-negative integer such that $\mseq-(k-\ell)\one\majleq\lseq+\ell\tseq$. 
Here $\one$ stands for the element $(1,\ldots,1)\in\mathbb{R}^{r+1}$.
\end{conjecture}


\bibliographystyle{bibsty-final-no-issn-isbn}
\addcontentsline{toc}{section}{References}
\ifthenelse{\equal{\finalized}{no}}{
\bibliography{abbrev,references}
}{
}
\end{document}